\documentclass[12pt,twoside]{amsart}

\usepackage{graphicx,mathrsfs,epstopdf}
\usepackage{amsmath,amsfonts,amssymb,amsthm}
\usepackage[colorlinks,citecolor=blue,urlcolor=blue]{hyperref}
\usepackage[dvipsnames]{xcolor}
\usepackage{algorithmic}
\usepackage[ruled,noend]{algorithm2e}
\usepackage{framed} 
\usepackage{blkarray}
\usepackage{pgfplots}
\usepackage{natbib}
\usepackage{enumerate}
\usepackage{subfigure}
\usepackage{tikz}
\usetikzlibrary{chains,matrix,scopes,fit,decorations,positioning,shapes,snakes,arrows,decorations.markings}
\newtheorem{thm}{Theorem}[section]

\newtheorem{prop}[thm]{Proposition}

\theoremstyle{definition}
\newtheorem{defn}[thm]{Definition}

\newtheorem{ex}[thm]{Example}

\theoremstyle{remark}
\newtheorem{rem}[thm]{Remark}

\DeclareMathOperator*{\argmax}{arg\,max}
\DeclareMathOperator{\rank}{rank}
\DeclareMathOperator{\Dir}{Dir}

\newcommand{\cC}{\mathcal{C}}

\newcommand{\E}{\mathbb{E}}

\newcommand{\cF}{\mathcal{F}}
\newcommand{\cG}{\mathcal{G}}

\newcommand{\N}{\mathbb{N}}

\renewcommand{\P}{\mathbb{P}}

\newcommand{\R}{\mathbb{R}}

\newcommand{\cX}{\mathcal{X}}

\newcommand{\bs}{\boldsymbol}

\newcommand{\cip}{\mbox{\,$\perp\!\!\!\perp$\,}}
\newcommand{\indep}{\cip}

 \newcommand\independent{\protect\mathpalette{\protect\independenT}{\perp}}
    \def\independenT#1#2{\mathrel{\rlap{$#1#2$}\mkern2mu{#1#2}}}
\newcommand{\notindependent}{\!\perp\!\!\!\!\not\perp\!}

\title[]{Tensors in algebraic statistics}

\author[]{Luis Sierra}
\address{}
\email{}
\thanks{}

\author[]{Marta Casanellas}
\address{}
\email{}
\thanks{}
\author[]{Piotr Zwiernik}
\address{}
\email{}
\thanks{}

\date{\today}                                            

\begin{document}
\begin{abstract}
Tensors are ubiquitous in statistics and data analysis. The central object that links data science to tensor theory and algebra is that of a model with latent variables. We provide an overview of tensor theory, with a particular emphasis on its applications in algebraic statistics. This high-level treatment is supported by numerous examples to illustrate key concepts. Additionally, an extensive literature review is included to guide readers toward more detailed studies on the subject.
\end{abstract}

\maketitle
\tableofcontents


\section{Introduction}


\noindent Tensors, which extend matrices to higher dimensions, were introduced in the late 19th century by authors such as William Rowan Hamilton, Woldemar Voigt and Gregorio Ricci-Curbastro, being made widely accessible by the work 
\cite{Riccilevicivita}. Although tensors were developed in differential calculus and geometry, nowadays tensors play a crucial role in statistics and data science as they offer a powerful mathematical framework for understanding complex phenomena, capturing intricate relationships, and extracting valuable insights from high-dimensional datasets (see the survey \cite{jietal19} for instance). In this overview paper, we explore the significance of tensors in statistics and data science, focusing on their mathematical foundations and connections to algebraic statistics. 

\noindent \textbf{Basic applications. }Tensors find extensive usage across various domains, including machine learning, neuroscience, and genomics. In machine learning, tensors facilitate the analysis of multidimensional data, enabling tasks such as image recognition, video summarization, and natural language processing. For example, in image recognition, tensors are used to represent high-resolution images, capturing spatial relationships and enabling the development of sophisticated deep learning models (see, e.g.,  \cite{Hartley2004,Bertolini24}). 

In neuroscience, tensors play a crucial role in unraveling brain connectivity patterns and understanding the functional organization of the brain (see for instance \cite{Cong15}). Brain networks, represented as higher-order tensors, capture the interactions among multiple brain regions, providing insights into cognitive processes, neurological disorders, and brain plasticity. Tensors enable the extraction of meaningful features from functional magnetic resonance imaging data, helping researchers identify brain regions associated with specific tasks or conditions, \cite{calamante2004}. 

Genomics is another domain where tensors offer a powerful framework for analyzing high-dimensional biological data. In multiomics studies, tensors are employed to integrate diverse molecular characteristics, such as gene expression, DNA methylation, and histone modifications (see for example \cite{Yang2018,Nafees20,Schreiber20}) . By modeling these complex relationships using tensors, researchers can identify molecular signatures associated with diseases, predict patient outcomes, and guide personalized treatment strategies.

\noindent \textbf{Tensor modelling. }The geometric complexity of tensors poses challenges in analysis and interpretation. To overcome this, researchers have developed methods that exploit simpler tensor structures. Many models
considered rely on more basic building blocks: diagonal tensors and rank-one tensors or low rank tensors. This simplification enhances interpretability and facilitates efficient computation and inference in high-dimensional settings.

As we will see, both the diagonal and the rank-one tensors arise naturally in the context of statistical independence. More generally, many of the probabilistic models that build on independence or conditional independence (e.g. Independence Component Analysis, Probabilistic Principal Component Analysis, Naive Bayes models) are based on finding a tensor decomposition into lower rank terms.

The field of algebraic statistics has emerged as a valuable tool for studying tensors within statistical models (bibliography on algebraic statistics includes the books \cite{Pistone2000, ASCB2005, drton2008lectures, sullivantbook, zwiernik2015semialgebraic}). By leveraging algebraic concepts, it provides insights into the properties and behavior of tensors in various statistical contexts. Algebraic tools have been successfully applied to address challenges related to the identifiability of latent variable models, parameter estimation, model fitting, design of experiments, and other fundamental statistical problems. 

\noindent \textbf{Latent variables and graphical models. } In data science, including \emph{latent} variables (or \emph{hidden}, unobserved) in the analysis is a popular way of improving expressibility of a model without dramatically increasing its computational complexity.  Despite their immense potential, the presence of latent variables introduces challenges in model estimation and interpretation. Addressing these challenges often requires sophisticated algorithms and statistical techniques. This formed an important motivation behind the development of the singular learning theory \citep{watanabe2009algebraic} and large part of recent development in algebraic statistics is centered around related problems.

By harnessing the algebraic structure of tensors, researchers have gained a deeper understanding of statistical models with latent variables and developed efficient estimation algorithms. Modelling dependence between random variables (both observed and unobserved) is usually done via a \emph{graphical model}. On a graphical model, the vertices of a graph represent random variables and its edges encode the conditional dependence relations. When this representation uses directed edges on an acyclic graph, the model is also known as a \emph{Bayesian network} or \emph{directed acyclic graph model}. The models for undirected graphs are sometimes called Markov random fields.

Latent trees form the most tractable family of Bayesian networks with unobserved variables, which can be used to model dependence structures when unobserved confounders are expected ~\cite[see for example][Section 2]{pearl2000}. However, there are several other reasons why latent tree models have become popular across sciences. First, latent tree models encompass a broader range of probability distributions compared to fully-observed tree models, while still leveraging their computational benefits, especially in high-dimensional scenarios. Additionally, the max-product algorithm enables efficient inference of unobserved states. Second,  trees can represent evolutionary processes and, in this context they are used in  phylogenetic analysis (see ~\cite{SteelPhylogeny} and Section \ref{sec:phylo}, in linguistics~\citep{ringe2002indo,shiers2014gaussian}, modeling language evolution, and in network tomography~\citep{castro2004network,eriksson2010toward} for inferring internet structures. Gaussian latent tree models (e.g., Brownian motion tree models) naturally capture correlations for network tomography due to correlations diminishing with tree distance; refer to Section~\ref{sec:gaussian}. Third, latent tree models capture hierarchical structure in complex datasets, and they are closely related to various hierarchical clustering methods \citep{bishop1998hierarchical,lawrence2004gaussian,zhangCluster,mourad2013survey}. This framework finds also application in computer vision as in~\cite{willsky2002multiresolution} and \cite{choi2010exploiting}.

Latent tree models generalize hidden Markov models (HMMs), which are defined on caterpillar trees~\cite{rabiner1989tutorial}. Hidden Markov Tree models (HMTMs) \cite{crouse1998wavelet} relax HMM's restrictions, admitting any unobserved tree structure. Applications span signal processing~\cite{choi2001multiscale,romberg2001bayesian}, biomedicine~\cite{makhijani2012accelerated,pfeiffer:CNE23820}, and linguistics~\cite{vzabokrtsky2009hidden}. This inclusiveness extends to phylogenetic models, Naive Bayes models, Brownian motion tree models, and one factor analysis model, leading to a unified framework for diverse model classes~\cite{wainwright2008graphical}.

An important example of latent variable model is Factor Analysis, which is commonly employed in the study of multivariate data to identify underlying latent factors that influence the observed variables. It finds applications in finance, psychology, and social sciences among others. A related model class is given by Principal Component Analysis (PCA) models, a widely used technique for dimensionality reduction and data compression. The probabilistic version, known as probabilistic PCA, introduces latent variables to account for the noise in the data, making it more flexible in real-world applications, see Section \ref{sec:ic}.

Another important class of graphical models with latent variables is the Restricted Boltzmann Machines (RBMs, see Section \ref{sec:RBM}). They serve as crucial components in deep learning, particularly in training deep belief networks and generative models. RBMs involve latent variables that capture hidden patterns in the data and can be extended to deep probabilistic graphical models with multiple layers of latent variables. They have had a transformative impact on fields like computer vision and speech recognition. From an algebraic standpoint, this model class has been thoroughly examined by Mont\'{u}far and collaborators \cite{montufar2015discrete,montufar2018restricted}.

\noindent \textbf{Summary. }
In this expository paper we explain the role of tensors in algebraic statistics by showcasing several examples. In Section~\ref{sec:defs} we present tensors that arise naturally in statistics and present the first results that give answers to statistical questions in terms of tensor properties. In Section~\ref{sec:GMs} we give a brief overview of graphical models and depict examples of both directed and undirected graphical models with either discrete or continuous variables by putting special emphasis in latent tree models. We devote Section \ref{sec:gaussian} to gaussian variables, which needs a special treatment. By building upon the previous examples, in Section~\ref{sec:tensorsAS} we illustrate the concepts of tensor rank and tensor decomposition and explore a method of moments. Finally, Section~\ref{sec:apps} is devoted to further applications of tensors in algebraic statistics, where we focus on phylogenetic analysis, restricted Boltzmann machines, and structural equation models.

\section{First examples and definitions}\label{sec:defs}

In this section we briefly introduce the main objects: tensors, symmetric tensors, diagonal tensors, rank-one tensors as well as their statistical counterparts. All concepts can be found in one of the books \cite{landsberg, zwiernik2015semialgebraic,bocci-chiantini,nonlinear}. 

A real $r_1\times r_2$ matrix is a two-dimensional array  $A=(a_{i_1,i_2})$ with $i_1\in \{0,\ldots, r_1-1\}$, $i_2\in \{0,\ldots,r_2-1\}$ and $a_{i_1,i_2}\in \R$. Given  $m\in \N$ and integers $r_1,\ldots,r_m\geq 1$, consider an $r_1\times \dots \times r_m$ \textit{array} of real numbers $A=[a_{i_1,\dots,i_m}]$ with $i_j\in \{0,\dots, r_j-1\}$ for $j=1,\ldots,m$. As we now explain, these arrays are naturally identified with \emph{tensors}. 

We start with the case of matrices. Let $\{e^j_{i}\}_{i=0,\dots,r_j-1}$ be the standard basis of $\R^{r_j}$, $j=1,\dots,m$. A matrix $A\in\R^{r_1\times r_2}$  can be identified with a  bilinear map from $\R^{r_1}\times \R^{r_2}$ to $\R$, $(x,y)\mapsto x^\top A y$ (which sends the standard basis pair $(e^1_i,e^2_j)$ to $a_{i,j}$). 
If $E^{i,j}$ is the matrix with zero entries except for a one in position $(i,j)$, its corresponding bilinear map is denoted as the tensor $e^1_i\otimes e_j^2$. In this way, a matrix $A$ as above corresponds to the tensor $\sum a_{i,j}e^1_i\otimes e^2_j$ in the vector space of tensors $\R^{r_1}\otimes \R^{r_2}$ (thought of as the space of bilinear forms on $\R^{r_1}\times \R^{r_2}$). 

Given two vectors $v^1=(v^1_0,\dots,v^1_{r_1-1})\in \R^{r_1}$ and $v^2=(v^2_0,\dots,v^2_{r_2-1})\in \R^{r_2}$, their \emph{tensor product} is defined as the tensor $v^1\otimes v^2=\sum_{i_1,i_2}v^1_{i_1}v^2_{i_2}\, e^1_{i_1}\otimes e^2_{i_2}  \in \R^{r_1}\otimes \R^{r_2}.$ For example, if $v^1=(2,0) \in \R^2$ and $v^2=(1,1)\in \R^2$, then $v^1\otimes v^2=2e^1_0\otimes e^2_0+2e^1_0\otimes e^2_1$ which is different from $v^2\otimes v^1=2e^1_0\otimes e^2_0+2e^1_1\otimes e^2_0$.

More generally, an $r_1\times \dots \times r_m$ tensor is often defined as a multilinear map $\R^{r_1}\times\dots \times \R^{r_m} \longrightarrow \R$: the tensor $e^1_{i_1}\otimes\dots \otimes e^{m}_{i_m}$ corresponds to the multilinear map 
\[
\begin{array}{rcl}
\R^{r_1}\times\dots \times \R^{r_m} & \longrightarrow &\R  \\(e^1_{j_1},\dots,e^{m}_{j_m})& \mapsto &  \begin{cases}
    1 & \text{ if }  (j_1,\dots,j_m)=(i_1,\dots,i_m),\\
    0, & \text{ otherwise}
\end{cases}  
\end{array}
\]
and these tensors form the standard basis of the vector space $\R^{r_1}\otimes \dots \otimes \R^{r_m}$ of $r_1\times \dots \times r_m$ tensors. In this sense, the tensor $t=\sum_{i_1,\dots i_m}a_{i_1,\dots,i_m}e^1_{i_1}\otimes\dots \otimes e^{m}_{i_m}$ in this space corresponds to the $r_1\times\dots \times r_m$ array $A=[a_{i_1,\dots,i_m}]_{i_j=0,\dots,r_j-1}.$ Coordinates of tensors will be refer to coordinates in the standard basis.

As above, given $m$ vectors $v^1=(v^1_0,\dots,v^1_{r_1-1})\in \R^{r_1}$, $\dots$, $v^m=(v^m_0,\dots,v^2_{r_m-1})\in \R^{r_m}$, their \emph{tensor product} is the tensor $v^1\otimes \dots \otimes v^m=\sum_{i_1,\dots,i_m}v^1_{i_1}v^2_{i_2}\dots v^m_{i_m}\,e^1_{i_1}\otimes\dots\otimes  e^m_{i_m}.$

We say that a tensor $T$ in $\otimes^r\R^n=\R^n \otimes \stackrel{r}{\ldots}\otimes \R^n$ is \emph{symmetric} if $T_{i_1\ldots i_r}=T_{j_1\cdots j_r}$ whenever $(j_1,\ldots, j_r)$ is a permutation of $(i_1,\ldots,i_r)$. The set of all symmetric tensors in 
$\otimes^r\R^n$
is denoted by $S^r(\R^n)$. The dimension of $S^r(\R^n)$ is $\binom{n+r-1}{r}$ and $T\in S^r(\R^n)$ can be codified by its entries $T_{i_1\cdots i_r}$ for $1\leq i_1\leq \ldots\leq i_r\leq n$.

\subsection{Discrete measures as tensors}\label{2:disc_meas}
Consider a collection of $m$ random variables, each with $r_i$ possible states, $X_i\in \cX_i=\{0,\ldots,r_{i}-1\}$ for $i=1,\ldots,m$. The vector $X=(X_1,\ldots,X_m)$ takes then values in $\cX=\cX_1\times \cdots \times \cX_m$. The probability distribution $p=(p(x))_{x\in \cX}$ of $X$ can be then identified with a point in 
$$
\R^{\cX}:=\R^{\cX_1}\otimes\cdots\otimes  \R^{\cX_m},
$$ 
where $\R^{\cX_i}$ is a copy of $\R^{r_i}$ with coordinates indexed by $\cX_i$. The coordinates of the tensor space $\R^\cX$ are  denoted by $p_x$, $x\in \cX$. For example, if $m=2$, $r_1=r_2=2$ then each $X_i$ is binary and the coordinates of $p$ in $\R^\cX$ are $p_{00},p_{01},p_{10},p_{11}$, $p=p_{00}e_0\otimes e_0+p_{01}e_0\otimes e_1+p_{10}e_1\otimes e_0+p_{11}e_1\otimes e_1$ (here we omit the superscripts  as both spaces coincide).

By definition, all probability distributions lie in the simplex $\Delta^\cX\subseteq \R^\cX$, where
 \begin{equation}\label{eq:probability-simplex}
\Delta^{\cX}\quad :=\quad \{p\in\R^{\cX}: \,\,p_{x}\geq 0, \,\sum_{x\in \cX} p_{x}=1\}.
\end{equation}
Any statistical model for $X$ is, by definition, a family of probability distributions and hence, a family of points in $\Delta^{\cX}$. This gives a basic identification of discrete statistical models with subsets of the tensor space $\R^\cX$. 

In geometric considerations concerning tensors, it is often more convenient to work in the projective space. Denote by ${\rm Proj}(V)$ the projectivization of a vector space $V$ and let $\P_\R^{\cX}:={\rm Proj}(\R^{\cX})$. By definition, this is the set of all points $\R^{\cX}\setminus \{0\}$ identified if they lie on the same line through the origin. Note that 
the probability simplex $\Delta^\cX$ is bijective with the nonnegative part of $\P_\R^{\cX}$. The bijection takes a point $q\in \P_{\R}^{\cX}$ with nonnegative coordinates and sends it to $p=(p_x)_{x\in \cX}\in \Delta^{\cX}$ with $p_x=q_x/\sum_{y\in \cX} q_y$. 

We note that the projective embedding is convenient because complex dependence models are often defined only up to the normalizing constant (see, e.g. the models defined by \eqref{eq:factorp}).

\subsection{Independence model and rank one tensors}\label{ex:independence}

    For a collection of discrete variables $X=(X_1,\cdots, X_m)\in \cX$, we say that the components of $X$ are independent (or $X$ satisifes the \emph{full idependence model}) if the tensor $p\in \Delta^\cX$ representing the distribution of $X$ can be written as a tensor product of the individual distributions $p^i\in \R^{\cX_i}$, $p=p^1\otimes \dots \otimes p^m$. In other words,
\begin{equation}\label{eq:fact0}
    p_{i_1i_2\cdots i_m}\;=\;p^1_{i_1}p^2_{i_2}\cdots p^m_{i_m}\qquad \mbox{for all }x=(i_1,\ldots,i_m)\in \cX.
\end{equation}


To have a geometric counterpart to this construction, we say that a tensor $p\in \R^\cX$ is a \emph{rank one tensor} (or a \emph{decomposable or pure} tensor) if it is the tensor product of an $m$-tuple of vectors $p^i\in \R^{\cX_i}$ for $i=1,\ldots,m$. The model of independence for the vector $X$ is contained in the set of rank one tensors. When $m=2$, rank one tensors correspond to rank one matrices. Indeed, 
$p=p^1\otimes p^2$ is the tensor $\sum_{i,j}p^1_{i}p^2_{j}\, e^1_{i}\otimes e^2_j$ which can be identified with the rank one matrix 
\[\begin{pmatrix}
    p^1_{0}\\ \vdots \\ p^1_{r_1-1} 
\end{pmatrix}
\raisebox{5mm}{$\begin{pmatrix}
    p^2_0 & \dots &p^2_{r_2-1}
\end{pmatrix}$}\]
via the correspondence introduced above.

In terms of projective algebraic geometry, the set of rank one tensors corresponds to the real Segre variety 
 \[{\rm Seg}(\P_{\R}^{\cX_1}\times \dots\times \P_{\R}^{\cX_m}) \,\subset \, \mathbb{P}_{\R}^{\cX}={\rm Proj}(\R^{\cX_1}\otimes \dots \otimes \R^{\cX_m}),
\]
see \cite[10.5.12]{bocci-chiantini}. For example ${\rm Seg}(\mathbb{P}_{\R}^1\times \mathbb{P}_{\R}^1)$ is the embedding of $\P_{\R}^1\times \P_{\R}^1$ in $\P_{\R}^3={\rm Proj}(\R^2\otimes\R^2)$ given by the map 
\[\begin{array}{rcl}
    \P_{\R}^1\times \P_{\R}^1 &\longrightarrow & \P_{\R}^3\\
    ([p^1_0:p^1_1],[p^2_0:p^2_1]) &\mapsto & [p^1_0p^2_0:p^1_0p^2_1:p^1_1p^2_0:p^1_1p^2_1] \,.
\end{array}\]
This variety is composed of those points $[p_{00}:p_{01}:p_{10}:p_{11}]$ in $ \P_{\R}^{3}$ that satisfy equation $p_{00}p_{11}-p_{01}p_{10}=0$ and corresponds to rank one tensors $p^1\otimes p^2$ or matrices
\[\begin{pmatrix}
    p_{00} & p_{01}\\
    p_{10} & p_{11}
\end{pmatrix}\] of rank one. Analogously, if we view $\R^{\cX_1}\otimes\R^{\cX_2}$ as matrices, ${\rm Seg}(\P_{\R}^{\cX_1}\otimes \P_{\R}^{\cX_2})$ is the set of $r_1\times r_2$ matrices of rank one (defined by the vanishing of all $2\times 2$ minors). Actually, it can be easily proven that the nonnegative part of ${\rm Seg}(\P_{\R}^{\cX_1}\times \dots\times \P_{\R}^{\cX_m})$ is isomorphic to the full independence model.

\subsection{Moment tensors - discrete case}\label{sec:discretemoms}

Consider discrete measures as discussed in Section~\ref{2:disc_meas}. For $u=(u_1,\ldots,u_m)\in \N^m$ and a vector $x=(x_1,\ldots,x_m)\in \cX$ define the monomial
$$
x^u\;:=\;x_1^{u_1}\cdots x_m^{u_m},
$$
where we use the convention that $0^0=1$. If $X$ has distribution $p$, then the corresponding \emph{moment} is 
$$
\mu_u\;=\;\E X^u\;=\;\sum_{x\in \cX}p_x x^u.
$$
This is called a \emph{moment of order} $k$ if $k=u_1+\dots+u_m$. In particular, $\mu_{0\cdots 0}=1$. 

\begin{ex}\label{eq:2binary}
    Consider the binary case $(X_1,X_2)\in \{0,1\}^2$. We have $\mu_{00}=p_{00}+p_{01}+p_{10}+p_{11}=1$, $\mu_{10}=p_{10}+p_{11}$, $\mu_{01}=p_{01}+p_{11}$, $\mu_{11}=p_{11}$.
    \end{ex}

\begin{prop}
The map $\R^\cX\to \R^{\cX}$ defined by 
$$\mu_{u}\;=\;\E X^u, \qquad \mbox{for all }u\in \cX,$$
it is a linear bijection. It maps the set $\{\sum_{x\in \cX} p_x=1\}$ to the set $\{\mu_{0\cdots 0}=1\}$.
\end{prop}

There are several other possible ways to encode a discrete probability distribution other than using moments. Some examples are cumulants. From the geometric perspective, they offer a (nonlinear) change of variables, which may be more suitable to study the geometry of a particular model; see for example \cite{sturmfels2013binary,zwiernik2012cumulants,ciliberto2016cremona}. For instance, in Example~\ref{eq:2binary} we could use the cumulant transformation $\kappa_{00}=0$, $\kappa_{10}=\mu_{10}$, $\kappa_{01}=\mu_{01}$, and $\kappa_{11}=\mu_{11}-\mu_{10}\mu_{01}$. The image of the independence model under this mapping is a linear subspace given by $\kappa_{11}=0$.

\section{Graphical models}\label{sec:GMs}

Many models in statistics and machine learning involve modelling dependence between the various vector components. Examples are multifold and include Factor Analysis models, Restricted Boltzmann machines, or Naive Bayes model. Graphical models serve as a natural setting for such a multivariate modelling task, as they aim to represent dependence relationships using graphs. In this context, the vertices of a graph represent random variables, and its edges encode the conditional dependence relations. According to the nature of the edges, we can distinguish two main types of graphical models: those with directed edges, and those with undirected edges, each suitable for a different tasks. 

For both types of graphical models, it is convenient to start with the notion of \emph{factorized distribution}. Let $X=(X_1,\ldots,X_m)$ be a random vector taking values in $\cX$ as above and let $\cF$ be a set of subsets of $[m]$. We may then consider the set of distributions over $\cX$ that satisfy
\begin{equation}\label{eq:factorp}
p(x)\;=\;\frac{1}{Z}\prod_{F\in \cF} \phi_F(x_F),\qquad \mbox{ for all }x\in \cX,    
\end{equation}
where $\phi_F$ are some nonnegative functions on $\R^{\cX_F}$ called potentials and $Z$ is the normalizing constant. 

In the discrete case, when $\cX$ is finite, we could represent each distribution in \eqref{eq:factorp} as a point in $\P^{\cX}_\R$ (ignoring the normalizing constant). The family of all such points admits a monomial parametrization $\prod_{F\in \cF} \phi_F(x_F)$, and so, it forms a toric variety in algebraic geometry terminology; see \cite{sullivantbook} for more details. 

A basic example is given in \eqref{eq:fact0}, where the distribution factorizes into the individual components. In this case the factors can be assumed to correspond to marginal probability distributions themselves but we will not require this to be true in general. 
\begin{ex}[Ising model]
The Ising model considers distributions on $\cX=\{0,1\}^m$ defined by
$$
p(x)\;=\;\frac{1}{Z}\prod_{i=1}^m e^{h_i x_i} \prod_{i<j} e^{J_{ij}x_i x_j}\;=\;\frac{1}{Z}e^{h^\top x+\tfrac12 x^\top J x},
$$
where $h=[h_i]\in \R^m$, $J=[J_{ij}]\in \R^{m\times m}$ is symmetric with zeros on the diagonal,  and $Z$ is the normalizing constant
$$
Z\;=\;\sum_{x\in \{0,1\}^n}e^{h^\top x+\tfrac12 x^\top J x}.
$$
It is therefore a factorized distribution. In many applications $J$ is supported on a graph $G$, in which case $J_{ij}=0$ whenever there is no edge between $i$ and $j$ in $G$. From the algebraic perspective, this model is known as the binary graph model and has been studied in \cite{develin2003markov}.
\end{ex}

If $\cX$ is finite as in Section~\ref{2:disc_meas}, it is clear that \eqref{eq:factorp} imposes some algebraic restrictions on the probability tensor $p\in \R^\cX$. We study various forms of these restrictions later in the paper. At this point it is perhaps only worth noting that such factorization properties are equivalent to tensor network models that have been widely considered in many different industrial applications. The paper \cite{robeva2019duality} explores this connection in more detail.

\subsection{Directed Acyclic Graphs}\label{bayesiannetwork}

Suppose $\cG = \left(V, E\right)$ is a directed graph. When the graph contains a directed edge $\left(a,b\right)\in E$, also denoted as $\{ a\rightarrow b\}$, we say $a$ is a \emph{parent} of $b$, and $b$ is a \emph{child} of $a$. A directed graph $\cG$ is acyclic (DAG) if it contains no directed cycles. For a vertex $v\in V$ in a DAG $\cG$  denote as $\mathbf{pa}(v)$, the set of \emph{parents} of $v$. By  $\mathbf{de}(v)$ denote the set of \emph{descendants} of $v$ as the vertices $w$ for which there exists a directed path $\{v\rightarrow w\}\subseteq \cG$, and conversely we may also define $\mathbf{nd}(v)$ as the \emph{nondescendants} $V\setminus \left(\left\{v\right\}\cup \mathbf{de}(v)\right)$. 

Given a DAG $\cG=(V,E)$ with  $V=\{1,\ldots,m\}$, associate a random variable $X_v$ at each node $v\in V$. 
The probabilistic model associated to the random vector $X=(X_1,\ldots, X_m)$ with values in $\cX=\prod_{i=1}^m \cX_i$ and the DAG $\cG$ contains all the distributions on $\cX$ that satisfy the following recursive factorization:
    \begin{equation}\label{bnfact}
        p({x}) = \prod_{i=1}^{n}{p\left(x_i|\mathbf{pa}(x_i)\right)},
    \end{equation}
where $p(x_i|x_C)$ denotes the conditional probability density function of $X_i$ given $X_C=x_C$, given by 
$$
p(x_i|\mathbf{pa}(x_i))\;=\;\frac{p(x_i,\mathbf{pa}(x_i))}{p(\mathbf{pa}(x_i))}
$$
for all $x\in \cX$ such that $p(\mathbf{pa}(x_i))>0$. 

\begin{rem}
For notational simplicity we use the same letter to denote all the conditional densities. 
To avoid confusion, a state denoted in lower case with the same name as a random variable (in upper case) refers to that observation on the random variable.
\end{rem}

In this article we focus on the situation when the random vector $X$ is exclusively either discrete or continuous, but mixed situations are also of interest \cite{lau96}. 

Given three random variables $X,Y,Z$ we say  that $X$ is independent of $Y$ given $Z$ (denotes $X\indep Y|Z$), if the conditional distribution satisfies $p(x,y|z)=p(x|z)p(y|z)$ for all $x,y,z$. Equivalently, $p(x|y,z)=p(x|z)$. The factorization in \eqref{bnfact} has an important interpretation in terms of conditional independence known as the directed \emph{local Markov property}:
    \begin{equation}\label{markovproperty}
        X_v \indep X_{\mathbf{nd}(v)\setminus \mathbf{pa}(v)} \: | \: X_{\mathbf{pa}(v)}, \quad v\in V.
    \end{equation}

\begin{figure}
    \begin{tabular}{ccc}
    (a)
      \begin{tikzpicture}[
    node distance=1cm,
    mynode/.style={draw, circle, fill, inner sep=1.5pt},
    >=stealth'
]

\node[mynode, label=above:$1$] (A1) {};
\node[mynode, label=above:$2$, right=of A1] (B1) {};
\node[mynode, label=above:$3$, right=of B1] (C1) {};
\draw[->, line width=1pt] (A1) -- (B1);
\draw[->, line width=1pt] (B1) -- (C1);
\end{tikzpicture}  &  \quad (b)\begin{tikzpicture}[
    node distance=1cm,
    mynode/.style={draw, circle, fill, inner sep=1.5pt},
    >=stealth'
]
\node[mynode, label=above:$1$] (A1) {};
\node[mynode, label=above:$2$, right=of A1] (B1) {};
\node[mynode, label=above:$3$, right=of B1] (C1) {};
\draw[->, line width=1pt] (C1) -- (B1);
\draw[->, line width=1pt] (B1) -- (A1);
\end{tikzpicture}  & \quad (c)\begin{tikzpicture}[
    node distance=1cm,
    mynode/.style={draw, circle, fill, inner sep=1.5pt},
    >=stealth'
]

\node[mynode, label=above:$1$] (A1) {};
\node[mynode, label=above:$2$, right=of A1] (B1) {};
\node[mynode, label=above:$3$, right=of B1] (C1) {};
\draw[->, line width=1pt] (B1) -- (A1);
\draw[->, line width=1pt] (B1) -- (C1);
\end{tikzpicture} \\[.4cm] (d)
      \begin{tikzpicture}[
    node distance=1cm,
    mynode/.style={draw, circle, fill, inner sep=1.5pt},
    >=stealth'
]
\node[mynode, label=above:$1$] (A1) {};
\node[mynode, label=above:$2$, right=of A1] (B1) {};
\node[mynode, label=above:$3$, right=of B1] (C1) {};
\draw[->, line width=1pt] (A1) -- (B1);
\draw[->, line width=1pt] (C1) -- (B1);
\end{tikzpicture}   & \quad (e) \begin{tikzpicture}[
    node distance=1cm,
    mynode/.style={draw, circle, fill, inner sep=1.5pt},
    >=stealth'
]

\node[mynode, label=above:$1$] (A1) {};
\node[mynode, label=above:$2$, right=of A1] (B1) {};
\node[mynode, label=above:$3$, right=of B1] (C1) {};
\draw[-, line width=1pt] (A1) -- (B1);
\draw[-, line width=1pt] (B1) -- (C1);
\end{tikzpicture} & \quad (f)\begin{tikzpicture}[
    node distance=1cm,
    mynode/.style={draw, circle, fill, inner sep=1.5pt},
    >=stealth'
]

\node[mynode, label=above:$1$] (A1) {};
\node[mynode, label=above:$2$, right=of A1] (B1) {};
\node[mynode, label=above:$3$, right=of B1] (C1) {};
\draw[-, line width=1pt] (A1) -- (B1);
\end{tikzpicture} 
    \end{tabular}
    
    \caption{Various graphs on three nodes. }
    \label{fig:3dags}
\end{figure}
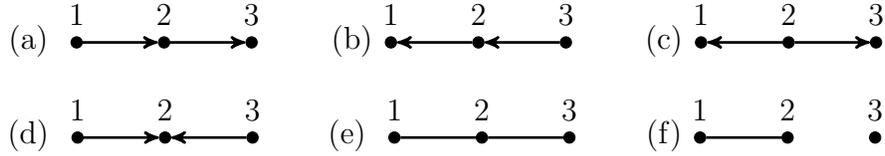

\begin{ex}
    Consider the structure of a \emph{Markov Chain}, in which random variables $X_1,\dots,X_m$ satisfy the property $p\left(x_k|x_{k-1}, \dots,x_1\right) = p\left(x_k|x_{k-1}\right)$ for all $k=2,\ldots,m$, so the future is only affected by the immediate present, and not by the past. Due to this conditional independence Markov property, a Markov chain can be interpreted as a probabilistic model on a line graph where each node represents a state (modelled by a random variable) $X_k$ with $X_{k-1}$ being its only parent, since $X_k$ is independent of the ``past'' given we know $X_{k-1}$. 
    If $m=3$, the associated graph is depicted in Figure~\ref{fig:3dags}(a).
\end{ex}

Analysis of DAG models leads to unexpected observations. For example, the graph (a) in Figure~\ref{fig:3dags} defines exactly the same model as the graphs (b) and (c). In particular, without additional assumptions, there is no one-to-one correspondence between DAG models and the corresponding directed graphs \cite{verma1990equivalence,chickering1995transformational}. 

\begin{ex}
    Consider a random vector consisting of three discrete random variables $X_1, X_2, X_3$ as components, where the dependence graph is given by Figure~\ref{fig:3dags}(d). In this model, the density factorizes as $p\left(x_1,x_2,x_3\right) = p\left(x_2|x_1,x_3\right)p\left(x_1\right)p\left(x_3\right)$. Observe that the marginal distribution of $(X_1,X_3)$ factors as $p\left(x_1,x_3\right) = p\left(x_1\right)p\left(x_3\right)$ and so naturally these variables are independent. By solving for the conditional distribution of $X_1,X_3$ given $X_2=x_2$ we may observe that this is no longer the case:
    \begin{equation*}
        p(x_1, x_3|x_2) = \frac{p\left(x_2|x_1,x_3\right)p\left(x_1\right)p\left(x_3\right)}{p(x_2)}.
    \end{equation*}
    Despite the fact that  $X_1 \independent X_3$, generally $X_1 \notindependent X_3 | X_2 = x_2$. Intuitively, this ``explaining away" phenomenon corresponds to the notion that despite two parents being independent, when looking at a child whatever is not explained by one of the parent must be explained by the other one, and so the independence is broken. Such a configuration is known as a \emph{collider} and is a well-known source of confusion in causal inference, sometimes referred to casually as Berkson's paradox.
\end{ex}

\begin{ex}\label{ex:3star}
Consider a simple DAG with four nodes and arrows $(0,1)$, $(0,2)$, and $(0,3)$, as in Figure~\ref{fig:singledag}. This corresponds to the factorization
 $$
 p(x)\;=\;p(x_0)p(x_1|x_0)p(x_2|x_0)p(x_3|x_0)\qquad\mbox{for all }x\in \cX.
 $$

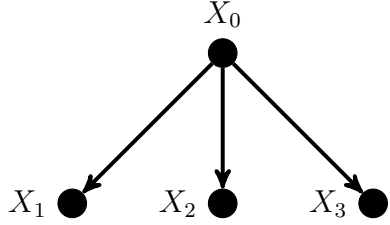
\begin{figure}
    \centering
    \begin{tikzpicture}[
    node distance=2.5cm,
    mynode/.style={draw, circle, fill, inner sep=0.4pt},
    >=stealth'
]
    \node[mynode, label=above:$X_0$] (1) at (0,0) {1};
    \node[mynode, label=left:$X_1$] (2) at (-2,-2) {2};
    \node[mynode, label=left:$X_2$] (3) at (0,-2) {3};
    \node[mynode, label=left:$X_3$] (4) at (2,-2) {4};
    
    \draw[->, line width=1.5pt] (1) -- (2);
    \draw[->, line width=1.5pt] (1) -- (3);
    \draw[->, line width=1.5pt] (1) -- (4);
    \end{tikzpicture}
    \caption{A graph on four nodes with a parent node and three descendants.}
    \label{fig:singledag}
\end{figure}
Suppose that $\cX=\{0,1\}^4$. Then each conditional distribution $p(x_i|x_0)$ for $i=1,2,3$ can be represented by a $2\times 2$ stochastic matrix (a nonnegative matrix with row sums all equal to $1$)
$$A^i\;=\;\begin{pmatrix}
     p(x_i=0|x_0=0) & p(x_i=1|x_0=0)\\
     p(x_i=0|x_0=1) & p(x_i=0|x_0=1)
\end{pmatrix}.$$
If we denote $\pi_x=p(X_0=x)$ and $a^i_{x,y}$ the $(x,y)$-entry of the matrix $A^i$, we see that the model is parametrized by
\begin{equation}\label{eq:param3star}     p_{x_0x_1x_2x_3}=\pi_{x_0}a^2_{x_0,x_1}a^3_{x_0,x_2}a^4_{x_0,x_3}.
\end{equation}
Note that the dimension of the space of 14 parameters $\pi_x$, $a^i_{x,y}$ is seven because we have to take into account that the sum of the rows of $A^i$ is one and $\pi_0+\pi_1=1$.
From \eqref{eq:param3star} it is easy to see that, if we fix the state at $X_0$ to be equal to $i$, the matrices
 \[\begin{pmatrix}
     p_{i000} & p_{i001}\\
     p_{i010} & p_{i011} \\
     p_{i100} & p_{i101}\\
     p_{i110} & p_{i111}
 \end{pmatrix} \qquad 
 \begin{pmatrix}
     p_{i000} & p_{i010}\\
     p_{i001} & p_{i011} \\
     p_{i100} & p_{i110}\\
     p_{i101} & p_{i111}
 \end{pmatrix} \qquad
\begin{pmatrix}
     p_{i000} & p_{i100}\\
     p_{i001} & p_{i101} \\
     p_{i010} & p_{i110}\\
     p_{i011} & p_{i111}
 \end{pmatrix}  
 \]
 have rank one.
From this, we easily get some of the defining equations of the model
$$
p_{ijkl}p_{ij'k'l'}=p_{ijkl'}p_{ij'k'l}=p_{ijk'l}p_{ij'kl'}=p_{ij'kl}p_{ijk'l'}\quad\mbox{for }i,j,j',k,k',l,l'\in \{0,1\}.
$$
In this way, conditional independence properties are translated into polynomial relationships among the entries of the probability vector $p$. In terms of algebraic geometry, we can think of the factorization \eqref{eq:param3star} as a monomial map
\[\varphi: \R^7 \longrightarrow \otimes^4\R^2\] which sends each set of free parameters to the distribution $p$. Both from the statistical and algebraic geometry perspective it is interesting to know a minimal generating set of polynomial equations that define the model (i.e that vanish precisely on the image of $\varphi$). In this case, it is known that the ideal of polynomials vanishing on ${\rm Im}(\varphi)$ is minimally generated by 18 quadratic equations of the above list, see \cite{garcia2005}.

This is the type of questions that lie at the core of algebraic statistics: for algebraic statistical models, use techniques of algebraic geometry and commutative algebra to find polynomial equations and semi-algebraic constraints that define the model. Knowing them might be useful for deciding whether the parameters are identifiable or not, for estimating parameters, for likelihood inference and for model fitting.
\end{ex}

\subsection{Undirected graphical models}

Let $\cG$ be an undirected graph and let $\cC$ denote the set of maximal cliques (complete subgraphs) of $\cG$.   An \emph{undirected graphical model} refers to the set of probability distributions
    \begin{equation}\label{ugmfact}
        f(x) = \frac{1}{Z}\prod_{C\in\mathcal{C}}{\psi_C(x_C)}, \qquad \textrm{where }Z = \int_{\mathcal{X}}\prod_{C\in\mathcal{C}}{\psi_C(x_C)}{\rm d}x,
    \end{equation}
    and $\psi_C:\cX_C\longrightarrow [0,\infty)$ are known as \emph{potential functions} due their interpretation in the context of statistical physics.

This definition may seem somewhat removed from the conditional independence notions described for directed graphs, but actually, a 1971 unpublished paper by John Hammersley and Peter Clifford showed that, under the additional restriction that $p(x)>0$ for all $x\in \cX$, the factorization \eqref{ugmfact} with respect to an undirected graph $\cG$ is equivalent to so called \emph{global Markov properties}: for any subsets $A,B,C\subset V$ we have $X_A\indep X_B|X_C$ as long as $C$ \textit{separates} $A, B$ in $\cG$, i.e. all paths between $A$ and $B$ pass through $C$.

Consider the discrete case as in Section~\ref{2:disc_meas}. Let $i_C:=(i_j)_{j\in C}\in\cX_C$ be a state of the vector $X_C$ for a clique $C\subseteq V$. Then, the potential function $\psi_C : \cX_C \longrightarrow [0,\infty)$ can be viewed as a vector $\theta^{(C)}$ with entries $\theta_{i_C}^{(C)}=\psi_C(i_C)\in [0,\infty)$. As in Example~\ref{ex:3star}, this gives as the parametrization
\begin{equation}
    p_{i_1 i_2 \cdots i_n} \; =\;  \frac{1}{Z\left(\theta\right)}\prod_{C\in\cC}{\theta_{i_C}^{(C)}}.
\end{equation}
In this way, we have again a parametric model with a semi-algebraic structure.

\subsection{Latent graphical models}\label{sec:latentGM}

Latent graphical models are designed to handle situations where certain variables $Z$ cannot be directly observed, yet they play a crucial role in explaining the relationships and patterns within the observed data $X$. In fact, using latent variables allows for very expressible models, meaning that they capture more of the complexity and variability of the data using a smaller number of covariates. For this reason, these \emph{latent} variables, also known as \emph{hidden} variables, are present in diverse domains and serve as essential components in modeling complex phenomena. That being said, they pose inferential difficulties, as we have previously pointed out.

\begin{ex}[Naive Bayes]\label{ex:naivebayes}
    Consider a classification problem over a dataset of discrete variables $X = (X_1, X_2, \dots, X_m)$ with another discrete variable $H$ encoding to which of $k$ classes $\mathcal{C} = \{c_1, c_2, \ldots, c_k\}$ an observation belongs. If $Z$ is unknown for new observations, such a class can be modeled as a latent variable. The Naive Bayes assumption is that the features are conditionally independent given the class, i.e., $X_i \independent X_j \:|\: Z$. The structure of the Naive Bayes model can be represented by a directed acyclic graph on the graph below:
    
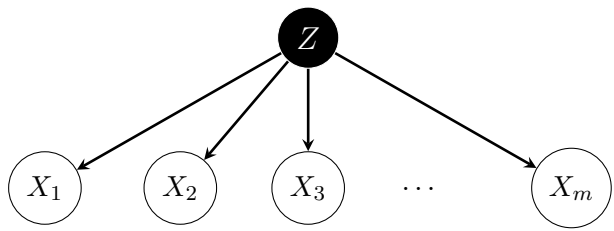
\vspace{5mm}
    \begin{center}
    \begin{tikzpicture}[>=stealth]
      \node[circle, draw, white, fill=black] (a1) at (0,0) { $Z$};
      \node[circle, draw] (a2) at (-3.5,-2) {\small $X_1$};
      \node[circle, draw] (a3) at (-1.7,-2) {\small $X_2$};
      \node[circle, draw] (a4) at (0,-2) {\small $X_3$};
      \node[] (a5) at (1.5,-2) {\small $\dots$}; 
      \node[circle, draw] (a6) at (3.5,-2) {\small $X_m$};
      
      \draw[->, line width=1pt] (a1) -- (a2);
      \draw[->, line width=1pt] (a1) -- (a3);
      \draw[->, line width=1pt] (a1) -- (a4);
      \draw[->, line width=1pt] (a1) -- (a6);
    \end{tikzpicture}
    \end{center}
   In order to predict the class $z \in \cC$ for a given observation of data $x = (x_1, \dots, x_m)$, we may use Bayes' rule to compute the posterior probability of the class:
    \begin{equation*}
        p(z | x) \propto p(z)\prod_{i=1}^{m}{p(x_i|z)},
    \end{equation*}
    and from the posteriors, the class $z \in \mathcal{C}$ is chosen by Bayes' decision rule which maximises said posterior probability:
     \begin{equation*}
        z = \argmax_{z\in\mathcal{C}}{\left\{p(z)\prod_{i=1}^{m}{p(x_i|z)}\right\}}.
    \end{equation*}
    This modeling approach is widely used in natural language processing and text classification tasks. Despite its simplistic assumptions, Naive Bayes often performs surprisingly well in practice, especially for tasks like spam detection, sentiment analysis, and topic categorization.
\end{ex}

\begin{ex}[Evolutionary model on a tripod tree]\label{ex:phylo3} 
As previously mentioned, latent tree models form an important example of latent variable models and such models have been a successful application of mathematics in biology. In Section~5, we consider Markov processes on phylogenetic trees, but for now we consider a simple case as a Naive Bayes model with $m=3$ where all random variables $X_1,X_2,X_3$ and $Z$ are discrete and take values in a state space $S=\{0,1,2,3\}$ representing the four nucleotides (adenine, cytosine, guanine, and thymine)
found in DNA sequences. Here the random variables $X_1,X_2,X_3$ correspond to observations of nucleotides at DNA sequences associated to biological living species (for example human, chimpanzee and gorilla) whereas the latent variable $Z$ codifies nucleotides in a DNA sequence of a common ancestor of these species. The local Markov property for this DAG gives
\[p(x_1,x_2,x_3,z)=p\left(z\right)p\left(x_1|z\right)p\left(x_2|z\right)p\left(x_3|z\right)\]
(cf. Example \ref{ex:3star}) and marginalizing over the latent variable we have
 \begin{equation}\label{eq:phylo3}
p\left(x_1,x_2,x_3\right)=\sum_{z\in S}p\left(z\right)p\left(x_1|z\right)p\left(x_2|z\right)p\left(x_3|z\right),
\end{equation}
  for any $(x_1,x_2,x_3) \in S^3.$ 
  These probabilities $p\left(x_1,x_2,x_3\right)$ stand for the probability of observing nucleotides $x_1,x_2,x_3$ at the leaves of the tree and can be estimated from observations of data at the living species, whereas $p(z)$ and $p(x_i|z)$ are (unknown) parameters of the model. In phylogenetics, estimating these parameters is relevant because they measure the amount of elapsed nucleotide substitutions from the common ancestor to each living species, which gives a notion of \emph{evolutionary distance} between an ancestral species and its descendants. In Section \ref{sec:tensorsAS} we explain methods to estimate these parameters from moments and in Section \ref{sec:apps} we explain these models in more detail. 
\end{ex}

\begin{ex}(Bag-of-words)\label{ex:bagofwords} 
   Consider now a restricted version of Naive Bayes where the $X_i$ are identically distributed conditioned on the latent variable $Z$. As a practical example, consider the bag-of-words, where a document composed of $m$ words has to be classified to one of some $k$ hidden \emph{topics}. We assume that the vocabulary has size $d$ 
    and each topic $c_h \in \mathcal{C}$ has its own probability vector $\mu_h\in \Delta^d$ for emitting words in the vocabulary: first a topic $c_h \in \mathcal{C}$ is chosen with probability $p(z=c_h)$, and then $m$ words $(x_1,\dots,x_m)$ are sampled independently with probability vector $\mu_h$.

    In this example, the ordering of the words in the document is not relevant, only their presence or absence matters, and therefore the words are assumed to be \emph{exchangeable} random variables, meaning their joint probability distribution $p(x_1, \dots, x_m)$ is invariant to the permutation of the indices and the corresponding tensor $p$ is a symmetric tensor. 
    
    Estimating parameters in a latent model is not an easy task. In Example~\ref{ex:mom_bag} we explain a specific method for the bag-of-words model, based on the method of moments. Another approach is to optimize the likelihood.
\end{ex}

\begin{rem}\label{em-latent}
As hinted in previous sections, a distribution $p(x, z|\, \theta)$ with latent variables where $z$ denotes a latent state, parameter inference can be challenging as $z$ is not directly observed. The canonical approach to estimating the parameters is that of looking to maximise the marginal likelihood over the observable variables 
$$
 {p(x|\theta)} \;=\; \sum_{z}{p(x,z|\,\theta)}, 
$$
where we can replace the sum by an integral if $Z$ also includes some continuous variables. 
Note that all the information known about the latent variables $Z$ is through their posterior distribution $p(\cdot|x,\theta)$. Therefore, since we do not have the likelihood for the complete dataset, it is instead estimated by taking the expected value of the posterior distribution (E-Step), and then this estimated likelihood is maximised with respect to the parameters $\theta$ (M-step); this procedure is known as the \emph{EM algorithm}, see \cite{dempster1977}.
\end{rem}

\section{The Gaussian case}

We devote this section to continuous random variables and, in particular, to normal multivariate variables. We start by generalizing moment tensors that were introduced above for the discrete case. 

\subsection{General moment tensors}
\label{sec:momenttensors}


Symmetric tensors typically arise as higher order derivatives of a smooth function. If $f:\R^n\to \R$ is such a function, we can arrange all order $r$ partial derivatives as a tensor $T\in S^r(\R^n)$ with coordinates
$$
T_{i_1\cdots i_r}\;=\;\frac{\partial^r }{\partial x_{i_1}\cdots \partial x_{i_r}} f(x)\qquad\mbox{for } 1\leq i_1,\ldots,i_r\leq n.
$$

Two types of functions and their partial derivatives will be particularly important in relation to statistics and machine learning. Consider the random vector $X=(X_1,\ldots,X_m)$, where now we allow the state spaces $\cX_i$ to be arbitrary subsets of $\R$. Let $M_X(\bs s) = \E e^{\bs s^\top X}$ and $K_X(\bs s)=\log \E e^{\bs s^\top X}$ denote the corresponding moment and cumulant generating functions, respectively. We consider the $r$-order $m\times \cdots\times m$ moment tensor $\mu_r(X)$, that is an $r$-dimensional table whose $(i_1,\ldots,i_r)$-th entry is 
$$
\mu_{i_1\cdots i_r}(X)\;=\;\E X_{i_1}\cdots X_{i_r}\;=\;\frac{\partial^{r} }{\partial s_{i_1}\cdots \partial s_{i_r}} M_X(\bs s)\Big|_{\bs s=0} ~.
$$ 
Similarly, the cumulant tensor $\kappa_r(X)$ is defined as having coordinates 
$$
\kappa_{i_1\cdots i_r}(X)\;=\;{\rm cum}(X_{i_1},\ldots,X_{i_r})\;=\;\frac{\partial^{r} }{\partial s_{i_1}\cdots \partial s_{i_r}} K_X(\bs s)\Big|_{\bs s=0} ~.
$$ 
The relationship between $\mu_r(X)$ and $\kappa_r(X)$ for arbitrary $r$ is cumbersome but very well understood \cite{speed1983cumulants,mccullagh2018tensor}. 
Directly by construction, $\mu_r(X)$ and $\kappa_{r}(X)$ are  symmetric tensors in $S^r(\R^m)$. 
\begin{rem}
Here we use a standard statistical notation for moment and cumulant tensors. In this convention we have, for example, that $\mu_1(X)=\E X_1$, $\mu_{13}(X)=\E X_1X_3$, $\kappa_{12}(X)=\E X_1X_2-\E X_1\E X_2$. This is in contrast with the special notation that we used in Section~\ref{sec:discretemoms} that is also routinely used in the algebraic study of discrete models.
\end{rem} 

For example, suppose that $X=(X_1,X_2)$ are two independent random variables. Then $\mu_{12}=\mu_{1}\mu_{2}$ and $\kappa_{12}=0$. 

In general, it may be not enough to study a statistical model only in the set of its low order moments/cumulants. An important exception is given by the multivariate Gaussian model (see next subsection), which is parameterized by the mean vector $\mu$ and the covariance matrix $\Sigma$.

\subsection{Multivariate Gaussian distribution}

The multivariate Gaussian distribution is the most important distribution in multivariate statistics and data science in general. It also plays an important role in many dimensionality reduction techniques for tensors, which we discuss in this paper. 

Denote by $S^2_+(\R^n)$ the set of all positive definite matrices. A $n$-variate Gaussian distribution with mean $\mu\in \R^n$ and covariance matrix $\Sigma\in S^2_+(\R^n)$ is a distribution on $\R^n$ with density
\begin{equation}\label{eq:normal}
f(x)\;=\;\det(2\pi \Sigma)^{-1/2} \exp(-\frac{1}{2}(x-\mu)^\top \Sigma^{-1}(x-\mu)).
\end{equation}
The moment and the cumulant generating functions of a Gaussian random vector $X$ are
$$
M_X(\bs s)\;=\;e^{\mu^\top \bs s+\tfrac12\bs s^\top \Sigma \bs s},\quad K_X(\bs s)\;=\; \mu^\top \bs s+\tfrac12\bs s^\top \Sigma \bs s.
$$
Note that the cumulant generating function is a quadratic polynomial and so all cumulants of order at least three are zero. By the classical result of J\'{o}zef Marcinkiewicz, the Gaussian distribution is the \emph{only} distribution for which the cumulant generating function is a polynomial \cite{marcinkiewicz1939propriete,lukacs1958some}. The fact that the cumulant generating function of a Gaussian distribution is a polynomial has important consequences in probability and statistics; see, e.g.,  \cite{janson1988normal}. The fact that no other distribution has this property is also important and it can explain the identifiability result of \cite{comon1994independent} that will come up in Section~\ref{sec:ic}.

One of the fundamental properties of the multivariate Gaussian distribution is that it is closed under taking margins and conditioning. More precisely, for an $n$-variate Gaussian vector $X$ with mean $\mu$ and covariance matrix $\Sigma$, consider an arbitrary split into two blocks $(X_A,X_B)$. Then the distribution of the sub-vector $X_A$ is Gaussian with mean $\mu_A=(\mu_i)_{i\in A}$ and covariance $\Sigma_{A,A}=[\Sigma_{ij}]_{i,j\in A}$. Moreover, the conditional distribution of $X_A$ given $X_B$ is Gaussian with the mean $\mu_A+\Sigma_{A,B}\Sigma_{B,B}^{-1}(X_B-\mu_B)$ and covariance $\Sigma_{A,A}+\Sigma_{A,B}\Sigma_{B,B}^{-1}\Sigma_{B,A}$. Thus working with marginal and conditional distributions of the multivariate Gaussian distribution corresponds to a purely algebraic operation on the parameters. 

\subsection{Gaussian graphical models}\label{sec:gaussian} The graphical model restrictions both for the directed and undirected graphs have a simple algebraic form when $X$ is assumed to follow a Gaussian distribution. 

Suppose that $\cG=(V,E)$ is an undirected graph with vertex set $V=\{1,\dots,m\}$ and $X=(X_1,\dots,X_m)$ is a multivariate Gaussian distribution. Then, by expanding \eqref{eq:normal} it is easy to see that the distribution of $X$ factorizes according to $\cG$ as in \eqref{ugmfact} precisely when the covariance matrix $\Sigma$ of $X$ satisfies
$$
(\Sigma^{-1})_{ij}\;=\;0\qquad\mbox{for all }(i,j)\notin E, i\neq j.
$$
Thus, parametrizing the model in terms of the inverse covariance matrix $K=\Sigma^{-1}$ allows us to define the model with simple linear restrictions. This observation motivated a line of work in algebraic statistics on understanding linear subspaces of symmetric matrices, e.g. \cite{uhler2010,sturmfels2020estimating,sturmfels2010multivariate, zwiernik2017maximum}.

Gaussian distributions that are defined by DAGs have also a simple description, which is equivalent to the underlying linear structural equation model (see Section~\ref{sec:ic}) by considering the system of linear stochastic equations
\begin{equation}\label{eq:linearsem} 
X_i\;=\;\sum_{j\in \textbf{pa}(X_i)}\lambda_{ij}X_j+\varepsilon_i,\qquad\mbox{for } i\in V,
\end{equation}
where $\varepsilon_i$ for $i\in V$ are assumed to be independent of each other and of the corresponding $\textbf{pa}(X_i)$. We assume that $\varepsilon_i$ have zero mean, are Gaussian, and the coefficients $\lambda_{ij}$ take arbitrary real values. In the matrix form, this can be written as $X=\Lambda X+\varepsilon$ with upper triangular $\Lambda$ (by reordering $X$). Since $\cG$ is acyclic, $I_m-\Lambda$ is invertible and we get $X=(I_m-\Lambda)^{-1}\varepsilon$. It follows that the inverse covariance matrix $K=\Sigma^{-1}$ satisfies
$$
K\;=\;(I_m-\Lambda)^\top \Omega^{-1} (I_m-\Lambda),
$$
where $\Omega$ is the diagonal matrix with variances of $\varepsilon_i$ on the diagonal. Thus, the model can be equivalently described by vanishing conditions on the Cholesky factorization of $K$. From the algebraic perspective, these models were studied in \cite{sullivant2008algebraic}.

\section{Tensor decomposition in algebraic statistics}\label{sec:tensorsAS}


Similar to matrix decomposition, tensor decomposition aims at representing a high-dimensional tensor as a sequence of elementary operations of simpler tensors. Tensor decompositions have been used in different fields, for example in computer vision, computational neuroscience, phylogenetics, or psychometrics. They are used to obtain blind source separation, independent component analysis and to provide estimators in latent variable models. We refer to \cite{jietal19} and \cite{Janzamin19} for good surveys on tensor decomposition and its applications in machine learning.

\subsection{Tensor rank}
It is an easy linear algebra exercise (using singular value decomposition for example) to see that an $n_1\times n_2$ matrix $M$ has rank $k$ if and only if it can be written as
\[M=u^1\otimes v^1+\dots+u^k\otimes v^k\] for some vectors $u^i\in \R^{n_1}$, $v^i\in \R^{n_2}$, 
and there does not exist such a decomposition of $M$ with a smaller number of summands (here $u^i\otimes v^i$ must  be viewed as the corresponding rank one matrix). 

A way to generalize the concept of matrix rank to tensors is the following.
\begin{defn}
 The \emph{rank} of a tensor $T\in \R^{r_1}\otimes\dots \otimes \R^{r_m}$ is the smallest integer $k$ such that $T$ can be written as 
 \begin{equation}\label{eq:rank}
     \sum_{j=1}^k u^{1,j}\otimes u^{2,j}\otimes\dots\otimes u^{m,j}
 \end{equation} for some $u^{i,j}\in \R^{r_i}$, $i\in [m]$, $j\in [k]$.
In other words, a tensor has rank $k$ if it can be written as a linear combination of $k$ rank one tensors and not as a sum of a smaller number of rank one tensors.

The \emph{symmetric rank} of a symmetric tensor $T\in S^r(\R^m)$ is defined as the smallest integer $k$ such that $T=\sum_{j=1}^k u^j\otimes u^j\otimes \dots \otimes u^j$ for some $u^j\in \R^m$, $j\in [k]$. 
\end{defn} 

We have already seen that tensors arising from the independence model of section~\ref{ex:independence} are rank one tensors. 
Below we provide examples of tensors of larger rank. In terms of algebraic geometry, as the algebraic variety $Seg\left(\P_{\R}^{X_{1}}\times \dots\times \P_{\R}^{\cX_m}\right)$ is the set of rank one tensors, the set of rank $k$ tensors lies in its $k$-th secant variety (and, actually, this secant variety is the Zarisky closure of this set of tensors).

\begin{ex}\label{ex:rknaive}
Tensors arising as the distribution of $X$ in the Naive Bayes model of Example~\ref{ex:naivebayes}, with a latent variable $Z$ taking values in a state space $\mathcal{C}$ of cardinal $k$ have rank at most $k$. Indeed, the Markov property gives the factorization of the joint distribution and, marginalizing over the latent variable $Z$, we get
\begin{equation*}
p(x_1,\dots,x_m)=\sum_{z\in \mathcal{C}} p(z)\prod_{i=1}^{m}{p(x_i|z)}
\end{equation*}
for any set of observations $x=(x_1,\dots,x_m)\in \cX$. Then, considering the vectors
\[u^{i,z}=\left(p(x_i=0|z),\dots, p(x_i=r_i-1|z)\right)\in \R^{\cX_i}\qquad\mbox{for }z\in \cC, i=1,\ldots,m,\]
we can write the distribution $p=\left(p{(x)}\right)_{x\in \cX}\in \R^{\cX}$ of $X=(X_1,\dots,X_m)$ as
\begin{equation}\label{eq:us}
p=\sum_{z\in \mathcal{C}} p(z) u^{1,z}\otimes\dots\otimes u^{m,z}.    
\end{equation}

Therefore, $p$ has tensor rank smaller than or equal to $k.$  In particular, the joint distribution $p$ at the leaves of the evolutionary model on the tripod tree of Example~\ref{ex:phylo3} is a tensor of rank smaller than or equal to four. 
\end{ex}

\begin{ex}
    Consider the bag of words model introduced in Example~\ref{ex:bagofwords}. 
    In this case, the vectors $u^{i,z}$ above coincide with $\mu_h$ if $z$ is the topic $c_h$, and do not depend on $i$. Thus, \eqref{eq:us} becomes
\[p=\sum_{h=1}^k p(z=c_h) \mu_{h}\otimes\dots\otimes \mu_{h}    
\]
and $p$ is a tensor of symmetric rank at most $k$.
\end{ex}

\subsection{Tensor decomposition}

A decomposition such as \eqref{eq:rank} has been introduced in different fields. It is called a \emph{tensor rank decomposition} or a \emph{Canonical Polyadic} (CP)  or \emph{Parallel Factor} (PARAFAC), see \cite{hitchcock27, Harshman1994}.

Such a decomposition is in general not unique. For example, one can rescale vectors $u_{i,j}$ or one can permute the summands in the decomposition. These are inherent ambiguities that can usually be sorted out depending on a particular application. One can restrict to special decompositions to obtain identifiability up to these ambiguities. As in the case of the spectral theorem for symmetric matrices, orthogonality is a common  requirement for decomposing symmetric tensors but not all symmetric tensors have an orthogonal decomposition. Uniqueness in the CP decomposition (up to permuting and scaling) can also be obtained for tensors of low rank, as we see below.

\begin{ex} Consider the Naive Bayes model of Example~\ref{ex:naivebayes} with $m=3$ and the corresponding probability tensor decomposition \eqref{eq:us}. Let $k$ be the cardinal of states of the latent variable. \cite{kruskal77} gave explicit conditions for the identifiability (up to permutation and scaling) of the summands in the tensor rank decomposition of the tensor given by the joint distribution. This explicit condition holds for \emph{generic} tensors in the model when $k$ is smaller than or equal to the cardinal $r_i$ of the state space of each observed variable. This contrasts sharply with the nonuniqueness of decompositions of rank $k$ matrices as sums of rank one matrices. \cite{AMR} showed how this basic result can be used in a powerful way to establish identifiability for a large class of latent variable models.

In the context of phylogenetic trees as in Example~\ref{ex:phylo3}, this result was rediscovered in \cite{chang1996full}. In this case the result says that the parameters of the model (the distribution of nucleotides at the latent variable and the conditional probabilities) can be uniquely recovered (up to permutation) from the joint distribution at the observed variables. Chang not only proved the identifiability of parameters but also gave an explicit construction for their recovery by using spectral decompositions of certain matrices obtained from the joint distribution. 
\end{ex}

\subsection{Nonnegative tensor rank}\label{nonnegativerank} The tensors that come from probability distributions, discrete or continuous, are a particular subset of the total class of tensors. For this reason, and especially for discrete variables, it is often useful to talk not about their rank as in \eqref{eq:rank}, but a variant called nonnegative tensor rank. A tensor $T$ has \emph{nonnegative rank at most} $k$ if it can be written as the sum of $k$ rank 1 tensors as before

\begin{equation}\label{tensorrank}
    T = \sum_{j=1}^{k}{u^{1,j}\otimes\cdots\otimes u^{m,j}},\quad u^{i,j}\in \R^{r_i}_{\ge0}, \: \: i\in\left[m\right],j\in\left[k\right],
\end{equation}
meaning that all of the vectors $u^{i,j}$ have nonnegative entries. 
The \emph{nonnegative rank} denoted $\rank_{+}$ is the smallest $k$ such that an expression like \eqref{tensorrank} holds. 
Directly by definition, $\rank T \le \rank_{+} T$ and, in general, the nonnegative tensor rank may be strictly larger. For more generalisations and variants of tensor ranks, see \cite{landsberg}.

\begin{rem}
    It is a simple exercise to check that the matrix given below has rank 3, but its nonnegative rank is equal to 4.
    \begin{equation*}
    \begin{pmatrix}
        1 & 1 & 0 & 0 \\
        0 & 1 & 1 & 0 \\
        0 & 0 & 1 & 1 \\
        1 & 0 & 0 & 1
    \end{pmatrix}
    \end{equation*}
    This is a simple instance of how the rank and nonnegative rank may differ.
\end{rem}

Determining whether the classical and nonnegative ranks coincide is an NP-hard problem. Nonnegative matrix factorization itself is an NP-hard problem, although polynomial performance is achievable via local search heuristics \cite{vavasis2018}. Thankfully, there are conditions under which such nonnegative rank decompositions are unique, and these are related to the identifiability of the underlying model \cite{qicommonlek2016}.

Nonnegative tensor decompositions have been studied in the context of PARAFAC by including nonnegativity constraints in chemometrics \cite{bro1997, bro1998}, and they can also be interpreted as a Naive Bayes decomposition of probability distributions
into conditionals as in \cite{garcia2005}.

Recalling the expression in \eqref{eq:fact0}, the independence model is an example of a nonnegative rank $\le 1$ tensor. For instance, the subset of $r_1 \times r_2$ matrices of $\rank_{+}\le 1$ is the nonnegative part of the Segre variety $\P^{\cX_1}\times \P^{\cX_2}$. Distributions on the independence model of two discrete random variables correspond to matrices in the intersection of $\Delta^{\mathcal{X}}$ with the set of $\rank_{+}\le 1$ matrices.

Mixture models make implicit use of nonnegative tensor rank. For a statistical model $\mathcal{P}\subset \Delta^{\cX}$, the $k$-th mixture model $\text{Mixt}^{k}\left(\mathcal{P}\right)$ is given by
\begin{equation*}
    \text{Mixt}^{k}\left(\mathcal{P}\right) = \left\{\sum_{i=1}^{k}{\pi_i p^{i}} \: : \: \pi \in \Delta^{k-1},\: p^{i}\in\mathcal{P}\right\}.
\end{equation*}
Such a model consists of all probability distributions of nonnegative rank less than or equal to $k$, and it is also a graphical model with hidden variables. A full treatment on this relationship can be found in \cite{drton2009}.

An important case is that of tensors of nonnegative rank $\le 2$, which are key in the study of Markov models on phylogenetic trees with binary states. The  Naive Bayes model in Example~\ref{ex:naivebayes} with binary observables and two latent classes can be represented by $2\times2\times\dots\times 2$ tensors of nonnegative rank at most $2$ which correspond to mixtures of 2 distributions. See \cite{allman2015tensors} for a characterization of such $\rank_{+} \le 2$ tensors and their geometry.

\begin{ex}[Topic Models] The model described in Example~\ref{ex:bagofwords} is a classifier based on document content, but it is not a generative model, and is also limited in its representational power. Topic models \cite{blei2012} are used when documents can exhibit themes from more than one topic, and classification would be too stringent of an approach. An example of such a model is that of Latent Dirichlet Allocation (LDA) \cite{lda2003}, which is a generative model that works as a Bayesian Network from Section~\ref{bayesiannetwork}, where one has a collection of $n$ documents and $k$ topics, with the property that elements are allowed to belong to various classes; a document can be about more than one topic. 

Topics are represented as probability distributions over words, and documents are represented in turn as distributions over latent topics. In fact, for each document, a scalar $\theta_{ij}\in\left[0,1\right]$ is given to each topic $j\in\left\{1,2,\dots,k\right\}$ according to how well the document fits into that topic ( that is, the probability of topic $j$ occurring in document $i$) as a mixture model where we assume that the total number of topics $k$ is known. The weights associated to the document topics follow a Dirichlet distribution $\Dir(\alpha)$, with probability density defined as
\begin{equation}
    p(\theta_i;\alpha) = f(\theta_{i1},\dots, \theta_{ik}; \alpha_1, \dots, \alpha_k) = \frac{1}{\mathbf{B}(\alpha)}\prod_{j=1}^{k}{\theta_{ij}^{\alpha_j-1}}
\end{equation}
where $\mathbf{B}$ is the Beta function, which acts as a normalizer. The previously mentioned weights $\theta_i = \left(\theta_{i1},\dots,\theta_{ik}\right)\in\Delta^{k-1}$ add up to 1, and as such are a discrete measure as in Section~\ref{2:disc_meas}. Then, we have that a word $i$ has probability $\varphi_{ij}$ of belonging to topic $j$, so similarly $\varphi_{ij}$ follows a distribution $\Dir(\beta)$.
The three sets of random variables representing documents, topics, and words, can be related by their conditional independence as a graphical model. By concatenating these three sets variables we naturally obtain a Bayesian network with the variables in the previously stated order as layers from parent to child.
\end{ex}

\subsection{Method of moments}

Learning the type of models discussed in Section~\ref{sec:latentGM} leads to many issues. The likelihood methods, which became a canonical choice for many simpler models, quickly become intractable in the presence of latent variables. With few good alternatives, the method of moments \cite{pearson1936method} is now routinely used in Gaussian mixtures \cite{amendola2016moment}, and latent tree models \cite{anandkumar2012method}, \cite{anandkumar14}, \cite{ruffini}. Since this method is intrinsically algebraic, we describe it more in detail in this manuscript. 

The idea of the method is simple. Suppose that the model parametrization induces an explicit parametrization of some moments $\mu=f(\theta)$ of the underlying distribution. A method of moments estimator then proceeds in two steps. First, a good sample estimator $\hat \mu$ of the population moments is found and then $\hat \theta$ is defined as a solution to $\hat\mu=f(\theta)$.

An obvious problem with this method is that, for a given $\hat\mu$, there is no guarantee that the equation $\hat\mu=f(\theta)$ has a solution. To some extent, this can be fixed by carefully choosing the moments for the analysis. But even if a solution exists, it may not satisfy some obvious restrictions. For this reason the method of moments is typically applied on a case by case basis and it requires a good understanding of the model class under consideration.

For instance, in the classic work \cite{pearson1894contributions}, a mixture of two Gaussian variables was used to model the body parts' sizes of a population of crabs. To do so, Pearson came up with a system of equations involving the first five moments of the distribution, in which he managed to eliminated the variables by obtaining a degree nine polynomial in a single variable (the product of the means). By solving the nonic, the parameters of the mixture distribution were obtained in this canonical example of the method of moments. For more details see \cite{amendola2018tapas}.

\begin{ex}\label{ex:mom_bag}(Method of moments for bag of words) 
Consider the bag of words in Example~\ref{ex:bagofwords} and assume that the number $m$ of words in the document is at least 3. 
\cite{anandkumar2012method} showed that using moments of order at most three, one can recover the parameters of the model, namely, the conditional probability vectors $\mu_h$ and the probabilities $w_h=p(z=c_h)$  for $h=1,\dots, k$. Here we explain this method of moments briefly following \cite[Appendix D]{anandkumar14}.

We can represent words in the vocabulary in correspondence with the standard basis $e_1,\dots,e_d \in \R^d$, and we consider $X$ so that the state space of $X_i$ is $\mathcal{D}=\{e_1,\dots,e_d\}$. 
Consider the variables $X_i$ as indicator vectors: 
\[X_i = e_j \iff \text{the $i$-th word in the document is $j$}.\]
Then the cross moments of these random vectors correspond to joint probabilities:
\[\E\left[X_1\otimes X_2\right]=\sum_{i,j}p(x_1=e_i,x_2=e_j)e_i\otimes e_j\]
and we get  
\[\E\left[X_1\otimes X_2\right]=\sum_{h=1}^k\sum_{i,j} p(z=c_h)p(x_1=e_i|z=c_h)p(x_1=e_j|z=c_h)e_i\otimes e_j.\] 
With the notation already introduced this is
\begin{equation}\label{eq:E2}
\E\left[X_1\otimes X_2\right]=\sum_{h=1}^k w_h \mu_h\otimes \mu_h.
\end{equation}
Similarly, we get
\begin{equation*}
\E\left[X_1\otimes X_2\otimes X_3\right]=\sum_{h=1}^k w_h \mu_h\otimes \mu_h\otimes \mu_h.
\end{equation*}
If we let $V=\left(\mu_1|\dots|\mu_k\right)$ and $D=diag(w_1,\dots,w_k)$, then considering $M_2=\E\left[x_1\otimes x_2\right]$ as a matrix we can express \eqref{eq:E2} as
\begin{equation}\label{eq:M2}
  M_2=VDV^\top.  
\end{equation}
This is not an eigendecomposition of $\E\left[X_1\otimes X_2\right]$ (except when $V$ is orthogonal), but we can use spectral techniques incorporating the third order moments. Indeed, consider $M_3=\E\left[X_1\otimes X_2\otimes X_3\right]$ and for any $\eta \in \R^d$ let $M_3\bullet\eta$ be the 2-way tensor 
\[(M_3\bullet\eta)_{i,j}=\sum_l (M_3)_{i,j,l}\, \eta_l . \] 
Then, as matrices, $M_3\bullet\eta=V D\, d(\eta) V^\top$ where $d(\eta)=diag(\mu_1^\top\eta,\dots,\mu_k^\top \eta),$ and  if $M=(M_3\bullet\eta) M_2^{-1},$
we have 
$M=V d(\eta) V^{-1}.$ Taking $\eta$ with distinct entries, the columns of $V$ are determined by the eigenvectors of $M$ (up to scaling and permutation). Once $V$ is recovered, $w_1\dots,w_k$ are easily recovered from expression \eqref{eq:M2}. 

There are more robust approaches for parameter estimation using moments, but these go beyond the scope of this survey. It is worth noting that the parameters estimated by a method of moments can also be used as the initial step of an EM algorithm (see Remark~\ref{em-latent}); in this case EM often converges to a local maximum after a single step; see \cite{Zhang16}.
\end{ex}

\section{Applications of algebraic statistics}\label{sec:apps}

\subsection{Phylogenetic analysis}\label{sec:phylo}

\begin{figure}			
				\tikzset{every picture/.style={line width=0.75pt}} 
				
				\begin{tikzpicture}[x=0.75pt,y=0.75pt,yscale=-0.7,xscale=0.7]
				
				\draw [line width=0.75]    (74.63,144.34) -- (103.37,185.55) ;
				\draw [shift={(104.72,187.48)}, rotate = 55.11] [color={rgb, 255:red, 0; green, 0; blue, 0 }  ][line width=0.75]      (0, 0) circle [x radius= 3.35, y radius= 3.35]   ;
				\draw [shift={(74.63,144.34)}, rotate = 55.11] [color={rgb, 255:red, 0; green, 0; blue, 0 }  ][fill={rgb, 255:red, 0; green, 0; blue, 0 }  ][line width=0.75]      (0, 0) circle [x radius= 3.35, y radius= 3.35]   ;
				\draw    (103.35,189.4) -- (74.63,229.76) ;
				\draw [shift={(74.63,229.76)}, rotate = 125.43] [color={rgb, 255:red, 0; green, 0; blue, 0 }  ][fill={rgb, 255:red, 0; green, 0; blue, 0 }  ][line width=0.75]      (0, 0) circle [x radius= 3.35, y radius= 3.35]   ;
				\draw [shift={(104.72,187.48)}, rotate = 125.43] [color={rgb, 255:red, 0; green, 0; blue, 0 }  ][line width=0.75]      (0, 0) circle [x radius= 3.35, y radius= 3.35]   ;
				\draw    (107.13,187.41) -- (155.94,187.48) ;
				\draw [shift={(158.29,187.48)}, rotate = 0.08] [color={rgb, 255:red, 0; green, 0; blue, 0 }  ][line width=0.75]      (0, 0) circle [x radius= 3.35, y radius= 3.35]   ;
				\draw [shift={(104.78,187.4)}, rotate = 0.08] [color={rgb, 255:red, 0; green, 0; blue, 0 }  ][line width=0.75]      (0, 0) circle [x radius= 3.35, y radius= 3.35]   ;
				\draw    (188.17,144.34) -- (159.62,185.55) ;
				\draw [shift={(158.29,187.48)}, rotate = 124.71] [color={rgb, 255:red, 0; green, 0; blue, 0 }  ][line width=0.75]      (0, 0) circle [x radius= 3.35, y radius= 3.35]   ;
				\draw [shift={(188.17,144.34)}, rotate = 124.71] [color={rgb, 255:red, 0; green, 0; blue, 0 }  ][fill={rgb, 255:red, 0; green, 0; blue, 0 }  ][line width=0.75]      (0, 0) circle [x radius= 3.35, y radius= 3.35]   ;
				\draw    (159.64,189.4) -- (188.17,229.76) ;
				\draw [shift={(188.17,229.76)}, rotate = 54.75] [color={rgb, 255:red, 0; green, 0; blue, 0 }  ][fill={rgb, 255:red, 0; green, 0; blue, 0 }  ][line width=0.75]      (0, 0) circle [x radius= 3.35, y radius= 3.35]   ;
				\draw [shift={(158.29,187.48)}, rotate = 54.75] [color={rgb, 255:red, 0; green, 0; blue, 0 }  ][line width=0.75]      (0, 0) circle [x radius= 3.35, y radius= 3.35]   ;
				
				\draw (4.7,129.13) node [anchor=north west][inner sep=0.75pt]  [font=\footnotesize]  {human};
				\draw (-18.2,231.28) node [anchor=north west][inner sep=0.75pt]  [font=\footnotesize]  {chimpanzee};
				\draw (193.35,129.13) node [anchor=north west][inner sep=0.75pt]  [font=\footnotesize] {gorilla}; 
				\draw (190.85,232.28) node [anchor=north west][inner sep=0.75pt]  [font=\footnotesize]  {orangutan};
				\draw (47,156.32) node [anchor=north west][inner sep=0.75pt]  [font=\footnotesize]  {$M^{e_1}$};
				\draw (47,197.32) node [anchor=north west][inner sep=0.75pt]  [font=\footnotesize]  {$M^{e_2}$};
				\draw (177,156.32) node [anchor=north west][inner sep=0.75pt]  [font=\footnotesize]  {$M^{e_3}$};
				\draw (177.5,197.32) node [anchor=north west][inner sep=0.75pt]  [font=\footnotesize]  {$M^{e_4}$};
				\draw (123.5,165.32) node [anchor=north west][inner sep=0.75pt]  [font=\footnotesize]  {$M^{e_5}$};
				\draw (105,192) node [anchor=north west][inner sep=0.75pt]  [font=\footnotesize]  {$r$};			
				\end{tikzpicture}
				\caption{\label{fig_phylotree} A phylogenetic tree with leaf set $L=\{\text{human, chimpanzee, gorilla, orangutan}\}$, root node $r$, and row stochastic matrices associated to a Markov process.}
			\end{figure}
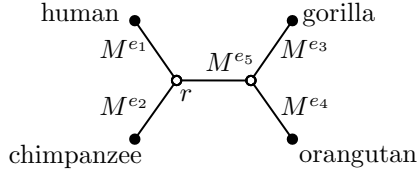

A \emph{phylogenetic tree} $T$ on a set $L$ of species (or other biological entities) is a tree with leaf set labelled by the elements of $L$. 
Figure~\ref{fig_phylotree} displays a phylogenetic tree with leaf set $L=\{\text{human, chimpanzee, gorilla, orangutan}\}$.

The main goal in phylogenetics is to reconstruct the evolutionary history of a set $L$ of living species from the information given by a collection of DNA molecules associated to hem. Due to the double strand symmetry, these molecules can be understood as words or sequences in the alphabet representing the four nucleotides. Under certain biological frameworks it can be assumed that each position in these sequences evolves independently of the other positions and in the same way.

If $T$ is a phylogenetic tree we model  the substitution of nucleotides on $T$ by first selecting an interior node $r$ to play the role of the root of the tree (this could represent the common ancestor to all species in $L$) and directing all edges out of it to view $T=(V,E)$ as a DAG. To each $v\in V$ we associate a random variable  $X_v$ taking values in $S=\{0,1,2,3\}$ (representing the four nucleotides $\{\texttt{A,C,G,T}\}$) and consider the local Markov property on the DAG. In this case, this property is equivalent to saying that each random variable is independent of its nondescendant variables given the observations at its immediate parent. It is well known that for the parameters to be identifiable one needs  that the tree has no nodes of degree 2. For $n=3$, this model is the one given in Example~\ref{ex:phylo3}.

To write the factorization of the joint distribution  \eqref{eq:fact0} under this Markov process one encodes all conditional probabilities in row stochastic matrices: we associate to each edge  $e=\{ u\rightarrow v\} \in E$ a $4\times 4$ matrix $M^e$ whose $x,y$ entry is $M^e_{x,y}=p(X_{v}=y| X_{u}=x)$, the conditional probability that a state $x$ at the ancestral node $u$ is being substituted by state $y$ at descendant node $v$. Then, 
the probability of observing states (nucleotides) $x_v$ at variables $X_v$, $v\in V$, is
\begin{equation}\label{eq:facttree}
    p(\{x_v\}_{v\in V})=p(x_r)\prod_{e=\{u\rightarrow u'\}\in E} M^{e}_{x_{u'},x_{u}}.
\end{equation}

In phylogenetics we only have observations for the random variables at the leaves of the tree and, as the interior nodes represent extinct species, the random variables at the interior nodes are latent. Thus, if the set of leaves $L$ is $[n]$, the probability $p_{i_1\dots i_n}$ of observing nucleotide $i_j$ at leaf $j$ for $j\in [n]$ can be obtained by marginalizing \eqref{eq:facttree} over the interior nodes:
\begin{equation}
   p_{i_1\dots i_n}=\sum_{\{x_v\}_{v\in V}\mid \, x_j=i_j, j\in [n]}   p(\{x_v\}_{v}).
\end{equation}
For each phylogenetic tree $T$ on $L=[n]$ we call $p^T\in \R^{\cX} =\R^4\otimes\stackrel{n}{\cdots}\otimes \R^4$ the tensor of the joint distribution 
$(p_{0\dots0},p_{0\dots1},\dots,p_{3\dots 3})$ 
obtained in this way. This model is known as the general Markov process on the phylogenetic tree.

If the phylogenetic tree that explains the evolutionary history is known and we have observations of DNA sequences on the set $L$, the parameters of the corresponding Markov process are usually estimated by a maximum likelihood approach. However, one of the main problems in phylogenetics is obtaining the tree $T$ that best fits the data: as the number of phylogenetic trees grows more that exponentially in the number of leaves, it is not possible to search the whole space of trees exhaustively. Here algebraic statistics plays an important role as we explain below (see also \cite{casanellas}).

In the late eighties, biologists Cavender, Felsenstein, and Lake realized that certain polynomial equations satisfied by the coordinates of $p^T$ can distinguish between different trees that can give rise to the distribution. That is, there are polynomial equations on the coordinates of a tensor $p\in \R^4\otimes\stackrel{n}{\cdots}\otimes \R^4$ that are satisfied if $p=p^T$ for some trees $T$ but not for others. One example of such equations comes from \emph{flattening} the tensor that we detail below.

Consider a  split of $L=[n]$ into two sets: a subset $A$ and its complement $B=L\setminus A$. As in section~\ref{ex:independence}, this split induces an isomorphism from the space of tensors to the space of matrices
\begin{equation}\label{eq:flatt}
\begin{array}{rcl}
  \R^{\cX} \cong \R^{\cX_A}\otimes \R^{\cX_B} & \longrightarrow & \mathcal{M}_{|\cX_A|\times |\cX_B|} (\R)\\
    p=(p_{x_1\dots x_n}) & \mapsto & {\rm flatt}_{A|B}(p)
\end{array},    
\end{equation}
where the $(x_A,x_B)$ entry of ${\rm flatt}_{A|B}(p)$ is the coordinate of $p$ corresponding to the probability of observing states $x_A=\{x_l\}_{l\in A}$ at the leaves in $A$ and $x_B=\{x_l\}_{l\in B}$ at the leaves in $B$. For example, in the phylogenetic tree of Figure~\ref{fig_phylotree}, if leaves human, chimp, gorilla and orangutan are denoted by $1,2,3,4$ respectively, we have
\[flatt_{12|34}(p)=\left(\begin{array}{ccccc}
p_{0000} & p_{0001} & p_{0002} &  \dots & p_{0033} \\
p_{0100} & p_{0101} & p_{0102} & \dots & p_{0133}\\
p_{0200} & p_{0201} & p_{0202} & \dots & p_{0233} \\
\vdots & \vdots & \vdots & \vdots & \vdots\\
p_{3300} & p_{3301} & p_{3302} & \dots & p_{3333}
\end{array}\right).\]

Assume that $p=p^T$ for some phylogenetic tree $T$ and consider a split $A|B$ of the set of leaves induced by removing an edge $e$ of $T$ (we call this an \emph{edge split}). Then ${\rm flatt}_{A|B}(p^T)$ can be considered as a joint probability table of the random variables $X_A=(X_l)_{l\in A}$ and $X_B=(X_l)_{l\in B}.$ By the Markov property on $T$ we know that $X_A \independent  X_B |X_u$ where $X_u$ is the hidden variable corresponding to one of the vertices of $e$. This is a Naive Bayes model with two observed random variables $X_A$, $X_B$ ($m=2$) and a latent one. Hence, $p^T$ is a tensor of rank at most four by Example~\ref{ex:rknaive}. In particular,  ${\rm flatt}_{A|B}(p^T)$ is a matrix of rank at most four because \eqref{eq:flatt} maps a sum of rank one tensors to a sum of rank one matrices. Thus, we have obtained the first claim of the following theorem:

\begin{thm}[\cite{allman2008}]
Let $T$ be a phylogenetic tree and let $A|B$ be a split of its set of leaves $L$. Let $p^T$ be a tensor distribution obtained from a Markov process on $T$. Then, if $A|B$ is an edge split on $T$,  ${\rm flatt}_{A|B}(p^T)$ has rank less than or equal to four. Moreover, if $A|B$ is not an edge split,  the rank of ${\rm flatt}_{A|B}(p^T)$ is larger than 4 (assuming that the parameters that generated $p^T$ satisfy certain genericity conditions).
\end{thm}

Therefore, the vanishing of the $5\times 5 $ minors of ${\rm flatt}_{A|B}(p)$ provide polynomial equations that are satisfied by distributions arising on trees that have $A|B$ as an edge split but are not satisfied for distributions on other trees. Instead of using these algebraic equations directly, it is more natural to evaluate the distance to rank four matrices by using singular values. This approach has been exploited to provide different successful methods of phylogenetic reconstruction, see for example  \cite{eriksson2005,chifmankubatko2014,fercas2016,casfergar2021}.

\subsection{Restricted Boltzmann Machines}\label{sec:RBM}

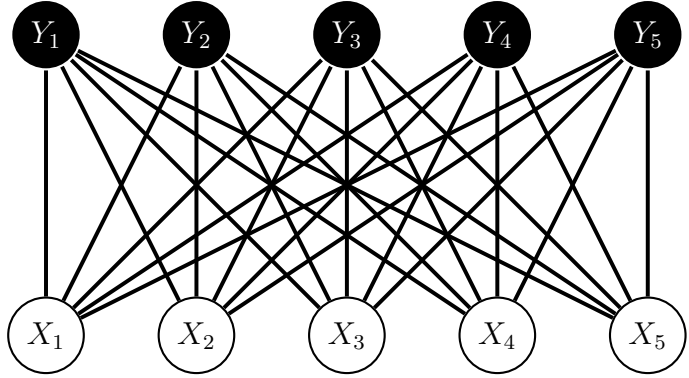
\begin{figure}
\begin{center}
\begin{tikzpicture}[thick]
    \foreach \i in {1,2,3,4,5}{
    \node[circle, draw, white, fill=black] (a\i) at (2*\i,0) {$Y_{\i}$};
    \node[circle, draw] (b\i) at (2*\i,-4) {$X_{\i}$};
}

\foreach \i in {1,2,3,4,5}{
\foreach \j in {1,2,3,4,5}{
    \draw[line width=1.5pt] (a\i) -- (b\j);
}
}

\end{tikzpicture}
\caption{\label{fig:RBM}A diagram showing a Restricted Boltzmann Machine as a graphical model with 5 latent and 5 observable variables.}
\end{center}
\end{figure}

Consider a set of binary random variables $V = \{X_1,\dots,X_n,Y_1,\dots,Y_m\}$ as the vertices of an undirected graphical model, where $X = \left(X_1,\dots,X_n\right)$ are observable variables and $Y = \left(Y_1,\dots, Y_m\right)$ are hidden variables. Such a graphical model is known as a Restricted Boltzmann Machines (RBM) and endowed a bipartite graph structure such where the only edges consist of one observed and one latent variable (see Figure \ref{fig:RBM}). The \emph{energy} of a given configuration of random variables with their parameters is defined as the function
\begin{equation}
    E(x, y; \theta) = -\sum_{i=1}^{n}{b_i x_i} -\sum_{j= 1}^{m}{c_j y_j} - \sum_{i=1}^n\sum_{j=1}^m W_{ij}x_i y_j,
\end{equation}
where $b\in \R^n$, $c\in \R^m$, and $W\in \R^{n\times m}$  are the parameters of the model; we write $\theta = \left(W, b, c\right)$. The notation for the parameters is influenced by the notation in the literature for perceptrons and neural networks. From the energy, one can define the probability of a realization of the random vector as 
\begin{equation}\label{prob_rbm}
    p(x,y;\theta) = \frac{1}{Z(\theta)}\exp(-E(x,y;\theta)), \quad \left(x,y\right)\in\left\{0,1\right\}^{n+m},
\end{equation}
with $Z(\theta)$ the normalising constant. As we may infer from the connectivity of the graph, the observable variables are independent with each other conditioned on the latent variables, and vice-versa. The marginal distribution of $X$ from equation \eqref{prob_rbm} is obtained by adding over the hidden variable states $y\in\left\{0,1\right\}^m$, and then via the appropriate manipulations we obtain the following product; see \cite{montufar2018restricted} for more details on the derivation and Boltzmann machines in general.

\begin{align}
    p(x;\theta) &= \frac{1}{Z(\theta)}\exp{\left(b^\top x\right)}\prod_{j=1}^{m}{\left(1+\exp{\left(w_{j}^\top x + c_j \right)}\right)} \nonumber \\ 
    &= \frac{1}{Z(\theta)}\prod_{j=1}^{m}{\left(\lambda_{j}\prod_{i=1}^{n}{{p_{ji}'\left(x_i\right)} + (1-\lambda_{j})\prod_{i=1}^{n}{{p_{ji}''\left(x_i\right)}}}\right)}\label{rbm_prod_str},
\end{align}
where $w_j$ is the $j$-th column of $W$ and  $p_{ji}',p_{ji}''$ are expressions which exclusively depend on the coordinate $x_i$. In this way we have an expression reminiscent of a product distribution $q(x_1,\dots,x_n)=\prod{q_i(x_i)}$. The link to tensor products is immediate. The expression in \eqref{rbm_prod_str} shows a product structure of nonnegative rank 2 tensors, which were highlighted in Section~\ref{nonnegativerank}. 
Indeed, in tensor notation \eqref{rbm_prod_str} can be written as
\begin{equation}
    p = \frac{1}{Z(\theta)}\prod_{j=1}^{m}{ \left( q_{j1}'\otimes \cdots\otimes q_{jn}' + q_{j1}'' \otimes \cdots\otimes q_{jn}''\right)}
\end{equation}
for some $q'_{ji}, q''_{ji}$.  As such, RBMs are an example of a product of $m$ factor tensors of nonnegative rank at most two. Such models have been explored and used for information processing in \cite{rbm1986smolensky} and in \cite{rbm1991haussler} as the precursor to two-layer neural networks, so their properties as function approximators are especially relevant in machine learning applications. Nonnegative rank two tensors have been described in terms of submodularity constraints in \cite{allman2015tensors}.

\begin{rem}
    The analysis of the distributions that arise from these rank two tensors is remarkably hard. However, for the simple case where $m=3$ and so the model is the sum of two $2\times2\times 2$ tensors, \cite{Seigal2017MixturesAP} showed that the model of nondegenerate distributions, that is those that lie on the interior of the space, are equivalent to a tree model with a single hidden variable with three states.
\end{rem}

\subsection{Independent Components and Structural Equations}\label{sec:ic}

Consider a model for a random vector $X$
\begin{equation}\label{eq:semc}
X\; = \;A\varepsilon,
\end{equation}
where $A\in \R^{m\times m}$ and $\varepsilon$ is a (latent) random vector with a simple dependence structure. Models of this form are extensively used in various applications. Depending on the context, we may impose additional restrictions on the matrix $A$ or the (latent) random vector $\varepsilon$. For example, in Blind Source Separation \cite[]{comon2010handbook}, it is commonly assumed that $A$ is invertible and that $\varepsilon$ has independent components. The following proposition will be important to this discussion:

\begin{prop}\label{prop:indep}
	The components of $X$ are independent if and only if $\kappa_r(X)$ is a diagonal tensor for every $r\geq 2$. 
\end{prop} 

Such a result has important applications in data science. In Independent Component Analysis, it is postulated that the observed random vector $X=(X_1,\ldots,X_m)$ can be written as a linear transformation $A\varepsilon$, where $\varepsilon=(\varepsilon_1,\ldots,\varepsilon_m)$ has independent components  and $A\in \R^{m\times m}$ is invertible. Under mild assumptions on the distribution of the (unobserved) vector $\varepsilon$, we can recover the matrix $A$ from observations of $X$ up to permuting its rows and flipping their signs \cite[]{comon1994independent}. In practice, this recovery is done with lower order cumulant tensors $\kappa_3(X), \kappa_4(X)$ using the fact that $\kappa_3(\varepsilon), \kappa_4(\varepsilon)$ are diagonal tensors. 

Since the vector $\varepsilon$ is not observed directly, the primary question is to what extent the matrix $A$ can be recovered solely from observations of $X$. Without loss of generality, we assume that $\varepsilon$ has zero mean and an identity covariance matrix. Under this assumption, the covariance of $X$ satisfies ${\rm cov}(X)=AA^\top$, which allows us to identify $A$ up to the orthogonal group action on the column space of $A$; specifically, $AA^\top=(AQ)(AQ)^\top$ for any orthogonal matrix $Q$. This identifiability result can be improved if $\varepsilon$ is non-Gaussian. Specifically, in the seminal paper \cite{comon1994independent}, it was shown that $A$ can be recovered up to permutation and scaling of the columns by $\pm 1$, provided that at most one component of $\varepsilon$ is Gaussian. Under mild generic conditions, $A$ can then be recovered from the third- or fourth-order moments of $X$. In practice, many methods utilize these lower-order moments.

To connect the model in \eqref{eq:semc} to our tensor discussion, let $\mu_r(X),\mu_r(\varepsilon)\in S^r(\R^m)$ denote the moment tensors of $X$ and $\varepsilon$, respectively, as discussed in Section~\ref{sec:momenttensors}. If $\varepsilon$ has independent components, then by Proposition~\ref{prop:indep}, $\mu_r(\varepsilon)$ is a diagonal tensor. Furthermore, by the multilinearity of moments, we have
$$\mu_r(X)\;=\;\mu_r(A\varepsilon)\;=\;A\bullet \mu_r(\varepsilon),$$
where the notation $A\bullet \mu_r(\varepsilon)\in S^r(\R^m)$ corresponds to the standard multilinear action of $A\in \R^{m\times m}$ on $S^r(\R^m)$:
$$
(A\bullet \mu_r(\varepsilon))_{i_1,\ldots,i_r}\;=\;\sum_{j_1,\ldots,j_r} A_{i_1j_1}\cdots A_{i_rj_r} (\mu_r(\varepsilon))_{j_1,\ldots,j_r}.
$$
In other words, \eqref{eq:semc} states that for every $r\geq 2$, the $r$-th moment tensor of $X$ is of the form $A\bullet D$ for some $A\in \R^{m\times m}$ and a diagonal tensor $D\in S^r(\R^m)$. The question then becomes: when does this condition, for a fixed $r$, allow us to identify $A$? Since the covariance matrix already allows us to identify $A$ up to the orthogonal group action on the column space, we can reduce this to the following problem. Suppose $A\bullet D=(AQ)\bullet D'$ for some invertible $A$ and orthogonal $Q$, with $D,D'$ being diagonal tensors. Equivalently, $Q\bullet D'$ is diagonal. Can we conclude that $Q$ must be a sign permutation matrix? For more detail on this approach, see \cite{mesters2022non}.
\begin{ex}
Consider the case when $r=3$, $m=2$. The condition that, for a diagonal tensor $D'$, the tensor $Q\bullet D'$ is diagonal, is given by the following system of two cubic equations in the entries of $Q$:
\begin{gather*}
(Q\bullet D')_{112}\;=\;Q_{11}^2Q_{21}D'_{111}+Q_{12}^2Q_{22}D'_{222}\;=\;0,\\
(Q\bullet D')_{122}\;=\;Q_{11}Q_{21}^2D'_{111}+Q_{12}Q_{22}^2D'_{222}\;=\;0.
\end{gather*}
In matrix form, this can be written as
$$
Q \cdot\begin{bmatrix}
Q_{11} & 0\\
0 & Q_{22}
\end{bmatrix}\cdot \begin{bmatrix}
Q_{21} & 0\\
0 & Q_{12}
\end{bmatrix}\cdot \begin{bmatrix}
D'_{111} \\
D'_{222}
\end{bmatrix}\;=\;\begin{bmatrix}
0 \\
0
\end{bmatrix}.
$$
Since $Q$ is orthogonal, each of the two diagonal matrices above is either identically zero or invertible. If either matrix is identically zero, then $Q$ must be a sign permutation matrix, and the equation clearly holds. If both are invertible, the equation cannot hold unless $D'_{111}=D'_{222}=0$, in which case $D'$ is the zero tensor. This shows that as long as $D'$ is a nonzero diagonal tensor, $Q\bullet D'$ is diagonal only if $Q$ is a sign permutation matrix.
\end{ex}

In most applications, we assume that $\varepsilon$ has independent components, but relaxations of these assumptions are also possible (e.g., the common variance model). For further motivation and basic results, see \cite{mesters2022non}. Additionally, the over-complete case, where $A$ is a rectangular matrix, is highly relevant in practice; see \cite{wang2024identifiability} for details.

Models of the form \eqref{eq:semc} have also been extensively studied in the context of causal analysis. In this setting, the matrix $A$ takes on a special structure. \cite{pearl2000} considered the following \emph{structural equations} model, defined on a directed graph:
$$
X_i = f_i(\mathbf{pa}(i), \varepsilon_i), \quad i = 1,\dots, m,
$$
which means that each variable in this system can be expressed as a function of its parent nodes and an error term $\varepsilon_i$. This system is nonlinear and nonparametric, but it generalizes the simpler linear version widely used in economics and the social sciences that was described in \eqref{eq:linearsem}. Compactly, we can write this as $X = \Lambda X + \varepsilon$, where $\Lambda_{ij} = 0$ unless $j\rightarrow i$ in $\cG$. By setting $A = (I_m - \Lambda)^{-1}$, we recover the model in \eqref{eq:semc}.

Given that $A$ has a special structure in this context, one might hope to identify it directly from the covariance matrix of $X$. For relevant results, see \cite{drton2011global,drton2016generic,foygel2012half}. The following example is adapted from \cite{drton2018algebraic}.
\begin{ex}
Does smoking during pregnancy affect the baby’s birth weight? To answer this question, suppose a study records the level of maternal smoking during pregnancy ($X_1$) and an infant’s birth weight ($X_2$). Assuming there is a causal effect of smoking on birth weight, we aim to quantify this effect. In practice, we center the data and often use a linear model:
$$X_1\;=\;\varepsilon_1,\qquad
X_2\;=\;\lambda_{12}X_1+\varepsilon_2,
$$
the goal is to infer the effect $\lambda_{12}$.  If $\varepsilon_1,\varepsilon_2$ are uncorrelated with mean zero, we can use the fact that ${\rm cov}(X_1,X_2)=\lambda_{12}{\rm var}(X_1)$ to recover $\lambda_{12}$ from the observed distribution.
\end{ex}

Another important generalization is to include latent variables in \eqref{eq:linearsem} that represent potential confounders. Such confounders may distort the direct causal relationships observed between the components of $X$; see \cite{barber2022half,tramontano2024parameter}. The following example is again taken from \cite{drton2018algebraic}.
\begin{ex}
The effect $\lambda_{12}$ in the previous example may not be direct and it may be distorted by some latent confounders. For instance, genetic or socio-economic factors might influence both smoking behavior and the baby's birth weight. In this case, inference based on the ratio ${\rm cov}(X_1, X_2)/{\rm var}(X_1)$ would no longer be valid. One solution to this issue is to introduce an ``instrumental variable'' that directly affects the level of smoking but not the birth weight and is not influenced by the latent variable. Following the literature, we could use the tax rate on cigarettes ($X_3$) in the state where each mother resides as such an instrumental variable. The linear equations defining the model with the latent variable $H$ are then given by:
\begin{gather*}
X_1\;=\;\lambda_{01}+\lambda_{31}X_3+\lambda_{H1}H+\varepsilon_1,\\
X_2\;=\;\lambda_{02}+\lambda_{12}X_1+\lambda_{H2}H+\varepsilon_2,\\
X_3\;=\;\lambda_{03}\;\;\quad\qquad\qquad\qquad+\varepsilon_3,\\
H\;=\;\lambda_{0H}\;\quad\qquad\qquad\qquad+\varepsilon_H,
\end{gather*} 
where the errors $\{\varepsilon_1,\varepsilon_2,\varepsilon_3, \varepsilon_H\}$ are pairwise independent. It follows that the equation ${\rm cov}(X_2,X_3) = \lambda_{12}\,{\rm cov}(X_1,X_3)$ holds in this model, and as long as ${\rm cov}(X_1,X_3) \neq 0$, we can infer properties of $\lambda_{12}$ based on the ratio ${\rm cov}(X_2,X_3)/{\rm cov}(X_1,X_3)$.
\end{ex}

Another problem that is highly relevant in causal analysis is  is whether the structure of the graph can be recovered from the observations of $X$. This again requires higher-order moment data, non-Gaussianity and is heavily based on the ICA identifiability result presented above \cite{shimizu2011directlingam}. 

\bibliographystyle{authordate1}

\end{document}